\documentclass{amsart}
\usepackage{amssymb}
\usepackage{latexsym}
\usepackage{amsmath}
\usepackage{euscript}

      \def\dC{{\mathbb C}}

   \def\dN{{\mathbb N}}   
      \def\dR{{\mathbb R}}

      \def\cC{{\mathcal C}}
\def\cD{{\mathcal D}}      
\def\cG{{\mathcal G}}   \def\cH{{\mathcal H}}   
   \def\cK{{\mathcal K}}   \def\cL{{\mathcal L}}
   \def\cN{{\mathcal N}}   \def\cO{{\mathcal O}}
      
      \def\cU{{\mathcal U}}

\def\h#1{{{\hat #1} }}

\def\bm\chi{\mbox{\boldmath$\chi$}}

\def\RE{{\rm Re\,}}

\def\ker{{\rm ker\,}}
\def\ran{{\rm ran\,}}

\def\dom{{\rm dom\,}}

\def\dim{{\rm dim\,}}

\let\xker=\ker \def\ker{{\xker\,}}

\def\max{{\text{\rm max}}}

\def\Real{{\text{\rm Re\;}}}
\def\Imag{{\text{\rm Im\;}}}

\def\ran{{\text{\rm ran\;}}}
\def\closp{{\text{\rm clsp\,}}}

\def\dom{{\text{\rm dom\;}}}

\def\leftik{{\bigl[\negthinspace\negthinspace\bigl[}}
\def\rightik{{\bigr]\negthinspace\negthinspace\bigr]}}
\def\senki{{\lbrack\negthinspace [\bot ]\negthinspace\rbrack}}
\def\senki+{{\lbrack\negthinspace [+] \negthinspace\rbrack}}
\def\senkit{{\lbrack\negthinspace [\cdot,\cdot ]\negthinspace\rbrack}}
\def\leftil{{\lbrack\negthinspace [}}
\def\rightil{{\rbrack\negthinspace ]}}

\newtheorem{theorem}{Theorem}[section]
\newtheorem{proposition}[theorem]{Proposition}
\newtheorem{corollary}[theorem]{Corollary}
\newtheorem{lemma}[theorem]{Lemma}

\theoremstyle{definition}

\newtheorem{definition}[theorem]{Definition}
\newtheorem{remark}[theorem]{Remark}

\numberwithin{equation}{section}

\newenvironment{proofspecial}%
{\begin{sloppypar}\noindent{\bf Proof of Theorem 3.1.}}%
{\hspace*{\fill}$\square$\end{sloppypar}\bigskip}

\begin{document}
\title[Nonlinear elliptic boundary value problems]
{A class of nonlinear elliptic boundary value problems}
\author[J. Behrndt]{Jussi Behrndt}

\address{Institut f\"ur Mathematik, MA 6-4 \\
Technische Universit\"at Berlin \\
Strasse des 17. Juni 136\\
10623 Berlin \\
Deutschland}
\email{behrndt@math.tu-berlin.de}

\maketitle

\begin{abstract}
In this paper second order elliptic boundary value problems on bounded domains $\Omega\subset\dR^n$ 
with boundary conditions on $\partial\Omega$ depending nonlinearly on the spectral parameter are investigated 
in an operator theoretic framework. For a general
class of locally meromorphic functions in the boundary condition a solution operator 
of the boundary value problem is constructed with the help of a linearization procedure. 
In the special case of rational Nevanlinna or Riesz-Herglotz functions on the boundary 
the solution operator is obtained in an explicit form in
the product Hilbert space $L^2(\Omega)\oplus (L^2(\partial\Omega))^m$, which is a natural generalization of 
known results
on $\lambda$-linear elliptic boundary value problems and $\lambda$-rational boundary value problems for
ordinary second order differential equations.
\end{abstract}


\section{Introduction}

Let $\Omega$ be a bounded domain in $\dR^n$, $n>1$, with smooth boundary $\partial\Omega$ and consider a uniformly
elliptic differential expression 
\begin{equation}\label{ellint}
\ell =-\sum_{j,k=1}^n\partial_j\,a_{jk}\,\partial_k \,+a
\end{equation}
on $\Omega$ with coefficients $a_{jk},a\in C^\infty(\overline\Omega)$ such that
$a_{jk}=\overline{a_{kj}}$ for all $j,k=1,\dots,n$ and $a$ is real-valued. The main objective of this paper is to solve 
the following eigenparameter dependent
boundary value problem: For a given function $g\in L^2(\Omega)$ and $\lambda$ in some open set $\cD\subset\dC$ 
find $f\in L^2(\Omega)$ such
that
\begin{equation}\label{bvpint}
(\ell -\lambda) f=g\qquad\text{and}\qquad \tau(\lambda)f\vert_{\partial\Omega}=\frac{\partial f_D}{\partial\nu_\ell}\bigl|_{\partial\Omega}
\end{equation}
holds. Here $\tau$ is assumed to be a meromorphic function on $\cD$ with values
in the space of bounded linear operators on $L^2(\partial\Omega)$, $\lambda$ is a point of holomorphy of $\tau$,
$f$ is a function in the maximal domain $\cD_\max=\{h\in L^2(\Omega):\ell h\in L^2(\Omega)\}$ 
and $f_D$ is the component of $f$ which lies in the domain of the Dirichlet operator. 

For the special case of a selfadjoint constant $\tau$ in the boundary condition in \eqref{bvpint} the boundary value problem
is uniquely solvable for all $\lambda$ which belong to the resolvent set of the selfadjoint partial differential operator
\begin{equation}\label{ttau}
T_\tau f=\ell f,\qquad \dom T_\tau=\left\{f\in \cD_\max: \tau f\vert_{\partial\Omega}=
\frac{\partial f_D}{\partial\nu_\ell}\bigl|_{\partial\Omega}\right\},
\end{equation}
in $L^2(\Omega)$ and the unique solution of \eqref{bvpint} is given by $f=(T_\tau-\lambda)^{-1}g$. Similarly, the nontrivial solutions of the
associated homogeneous problem, i.e., $g=0$ in \eqref{bvpint}, are given by the eigenvectors corresponding to the 
(real) eigenvalues $\lambda$ of $T_\tau$. 

Elliptic problems with $\lambda$-linear boundary conditions were already considered by J.~Ercolano and M.~Schechter in \cite{ES65-1,ES65-2}
and a solution operator $\widetilde A$ in the larger space $L^2(\Omega)\oplus L^2(\partial\Omega)$ was constructed
and its spectral properties were studied. Again the resolvent of $\widetilde A$, or, more precisely, the compression of 
the resolvent onto the basic space $L^2(\Omega)$, 
\begin{equation*}
f=P_{L^2(\Omega)}(\widetilde A-\lambda)^{-1}\!\upharpoonright_{L^2(\Omega)} g ,
\end{equation*}
yields the unique solution $f$ of \eqref{bvpint}, and the eigenvalues and the (components in $L^2(\Omega)$ of the) eigenvectors 
of $\widetilde A$ are the nontrivial solutions of the homogeneous problem. 
We emphasize that the solution operator $\widetilde A$ in the $\lambda$-linear case is selfadjoint with respect to
the Hilbert scalar product in $L^2(\Omega)\oplus L^2(\partial\Omega)$ if $\tau(\lambda)=\lambda$ and 
selfadjoint with respect to an indefinite (Krein space) inner product if $\tau(\lambda)=-\lambda$. 
The spectral properties of selfadjoint operators in Krein spaces differ essentially from the spectral properties of
selfadjoint operators in Hilbert spaces and this affects the solvability of \eqref{bvpint}. E.g., 
if $\tau(\lambda)=-\lambda$ in \eqref{bvpint}, then
the solution operator $\widetilde A$ and the homogeneous boundary value problem may have non-real 
eigenvalues, see \cite{BHLN01}. 

The main objective of this paper is to go far beyond the $\lambda$-linear case and to investigate the
solvability of the boundary value problem \eqref{bvpint} for a large class of 
operator-valued functions in the boundary condition. Here it will be assumed that $\tau$ is
a meromorphic function on some simply connected open set $\cD\subset\dC^+$ with values in the space 
$\cL(L^2(\partial\Omega))$ of bounded linear operators on $L^2(\partial\Omega)$ and that $\tau$ admits a 
minimal representation
\begin{equation}\label{oprep}
\tau(\lambda)=\Real\tau(\lambda_0)+\gamma^+\bigl((\lambda-\Real\lambda_0)+(\lambda-\lambda_0)
(\lambda-\bar\lambda_0)(A_0-\lambda)^{-1}\bigr)\gamma
\end{equation}
with the help of the resolvent of a selfadjoint operator or relation $A_0$ in a Krein or Hilbert space 
$\cH$ and a mapping 
$\gamma\in\cL(L^2(\partial\Omega),\cH)$. We mention that, e.g., locally holomorphic functions, Nevanlinna and generalized
Nevanlinna functions, and so-called definitizable and locally definitizable functions can be represented in the form 
\eqref{oprep}, see \cite{A84,DLS87,HSW98,J92,J00,J05,KL77,LT77}.

For the construction of a solution operator $\widetilde A$ of the boundary value problem \eqref{bvpint} 
we make use of the notion
of (generalized) boundary triples, and associated Weyl or $M$-functions, a convenient and useful tool for the spectral analysis
of the selfadjoint extensions of an arbitrary symmetric operator with equal deficiency indices, see, e.g., \cite{BGP08,D95,D99,DM91,DM95,GG91}. 
Boundary triplets for the maximal operator $T_\max f=\ell f$,
$f\in\cD_\max$, 
generated by the elliptic differential expression in $L^2(\Omega)$ were used (also in the non-symmetric case) in 
\cite{BGW09,GM08,GBook08} and appear in a slightly different form already in the fundamental paper \cite{G68} of G.~Grubb.
One of the main ingredients in the construction of a solution operator $\widetilde A$ of \eqref{bvpint} is to realize 
the function $\tau$
in the boundary condition as the Weyl function corresponding to some boundary triple, cf. \cite{B08,BL07,DHMS00}
and \cite{ACD03,BS09,CDR01,DLS87-1,DLS88,DLS93,DL96,LM91} for other approaches.
So far this is possible only under rather restrictive assumptions on the function $\tau$, e.g.,
in the special case of an $\cL(L^2(\partial\Omega))$-valued Nevanlinna function one has to assume that 
$\Imag\tau(\lambda)$ is boundedly invertible, see \cite{DM91,LT77}, or one has to apply the concept
of boundary relations and Weyl families from \cite{DHMS06,DHMS07}. 
Therefore, in order to treat the problem \eqref{bvpint} in a general setting, we extend the existing results on realizations of operator functions as Weyl functions in Section~\ref{sec2}.
Here a new method is proposed in which 
an arbitrary operator function $\tau$ of the form \eqref{oprep} can be realized 
as the Weyl function corresponding to a generalized boundary
triplet associated to a restriction of the selfadjoint operator or relation $A_0$. The idea is based on a
decomposition of $\tau$ in a constant part and a ``smaller'' part which satisfies a special strictness condition, 
see Definition~\ref{strictdef} and \cite{PAMS} for the special case of matrix Nevanlinna functions. 
Although the realization obtained in Theorem~\ref{mainthm} is in general not minimal it turns out that the connections between the 
solvability of the boundary value problem \eqref{bvpint} and the spectral properties of the solution operator $\widetilde A$ 
are not affected at all. 

The heart of the paper is Section~\ref{pde}, where the eigenvalue dependent boundary value problem \eqref{bvpint}
is discussed. After recalling some basic properties on elliptic operators associated to \eqref{ellint} 
and a corresponding ordinary boundary triple for $T_\max$ in Section~\ref{pde1} 
we construct a solution operator $\widetilde A$ of the elliptic boundary value problem \eqref{bvpint} in a larger 
Krein or Hilbert space $L^2(\Omega)\times\cK$ with the help of the realization result from Section~\ref{sec2}.  
The unique solution $f\in L^2(\Omega)$ of \eqref{bvpint} and the compression of the resolvent of $\widetilde A$
onto the basic space $L^2(\Omega)$ are then expressed in the form
\begin{displaymath}
f=P_{L^2(\Omega)}(\widetilde A-\lambda)^{-1}\!\upharpoonright_{L^2(\Omega)}g=(T_D-\lambda)^{-1}g-
\gamma(\lambda)\bigl(M(\lambda)+\tau(\lambda)\bigr)^{-1}\gamma(\bar\lambda)^*g,
\end{displaymath}
where $T_D$ is the Dirichlet operator asssociated to $\ell$ in $L^2(\Omega)$, $M$ denotes the Weyl or $M$-function
corresponding to an ordinary boundary triple for $T_\max$ and $\gamma(\cdot)$ is the associated $\gamma$-field, 
cf. Proposition~\ref{elbt}. We point out that for a constant selfadjoint boundary condition $\tau$ the solution operator $\widetilde A$
coincides with $T_\tau$ in \eqref{ttau} and the above formula reduces to the well-known Krein formula
for canonical selfadjoint extensions in $L^2(\Omega)$ of the minimal operator associated to $\ell$, cf. \cite{BL07,BGW09,GM08-1,GM08-2,GM08,G08,P08,PR09,Post07,R07}.
The proof of our main result Theorem~\ref{rwpthm} is based on a coupling technique of 
ordinary and generalized boundary triples which differs from the methods
applied in earlier papers. 

We illustrate our general approach in Section~\ref{tauex} in an example 
where $\tau$ is chosen to be 
a rational $\cL(L^2(\partial\Omega))$-valued Nevanlinna (or Riesz-Herglotz) function of the form
\begin{equation}\label{tauratint}
\tau(\lambda)=\alpha_1+\lambda\beta_1+\sum_{i=2}^m\beta_i^{1/2}(\alpha_i-\lambda)^{-1}\beta_i^{1/2}
\qquad\lambda\in\bigcap_{i=2}^m\rho(\alpha_i).
\end{equation}
Here $\alpha_i,\beta_i$ are bounded selfadjoint operators on $L^2(\partial\Omega)$ and $\beta_i\geq 0$.
In this special case the solution operator from Theorem~\ref{rwpthm} acts in the product space $L^2(\Omega)\oplus
(L^2(\partial\Omega))^m$ and can be constructed in a more explicit form, cf. Theorem~\ref{ratrwp} and Corollary~\ref{lincor}
for the $\lambda$-linear problem. We point out that 
an analogous selfadjoint solution operator in $L^2(I)\oplus\dC^m$ 
of a Sturm-Liouville problem on a bounded interval $I\subset\dR$ with a scalar variant of
\eqref{tauratint} in the boundary condition was constructed in \cite{B03}.

The paper is organized as follows. In Section~\ref{sec1} we give a brief introduction into the theory
of ordinary boundary triples and generalized boundary triples associated to symmetric operators and
relations in Krein spaces. The corresponding $\gamma$-field and Weyl function are defined 
and some of their basic properties are recalled. In Section~\ref{sec2} it is shown 
how an arbitrary operator function $\tau$
of the form \eqref{oprep} can be interpreted as the Weyl function of some generalized boundary triple and
some special classes of operator functions are discussed in Section~\ref{sec2.3}. Section~\ref{pde}
treats the elliptic boundary value problem \eqref{bvpint}, in particular, a solution
operator $\widetilde A$ is constructed, it is shown that the compressed resolvent of $\widetilde A$ onto the
basic space $L^2(\Omega)$ yields the unique nontrivial solution 
of the inhomogeneous problem \eqref{bvpint} and that the eigenvalues and eigenvectors of $\widetilde A$ 
solve the homogenous boundary value problem.

\section{Generalized boundary triples and Weyl functions of symmetric
relations in Krein spaces}\label{sec1}

Let $(\cH,[\cdot,\cdot])$ be a Krein space and let $J$ be a corresponding fundamental symmetry.
We study linear relations in $\cH$, that is, linear
subspaces of $\cH\times\cH$. The elements in a linear relation will be denoted by $\hat f=\{f,f^\prime\}$, 
$f,f^\prime\in \cH$.
For the set of all closed linear relations
in $\cH$ we write $\widetilde\cC(\cH)$. Linear operators in
$\cH$ are viewed as linear relations via their graphs. 
The linear space of bounded linear
operators defined on a Krein space $\cH$ with values in a Krein
space $\cK$ is denoted by $\cL(\cH,\cK)$. If $\cH=\cK$ we simply
write $\cL(\cH)$. We refer the reader to \cite{AI89,B74,DS87-1,DS87-2} for more details on Krein spaces and linear operators and
relations acting therein.

We equip $\cH\times\cH$ with the Krein space inner product $\senkit$ defined by
\begin{equation}\label{ind}
\leftik\hat f,\hat g\rightik:=i\bigl([f,g^\prime]-[f^\prime,g]\bigr),\qquad\hat f=\{f,f^\prime\},
\,\hat g=\{g,g^\prime\}\in\cH\times\cH.
\end{equation}
Then $\bigl(\begin{smallmatrix} 0 & -iJ\\ iJ & 0 \end{smallmatrix}\bigr)\in\cL(\cH^2)$ is a corresponding 
fundamental symmetry. Observe that also in the special case when $(\cH,[\cdot,\cdot])$ is a
Hilbert space, $\senkit$ is an indefinite metric. In the following we shall often use at the
same time inner products $\senkit$ arising from different Krein and Hilbert spaces
as in \eqref{ind}. Then we shall indicate these forms by subscripts, for example,
$\senkit_{\cH^2}$, $\senkit_{\cG^2}$.

For a linear relation $A$ in the Krein space $\cH$ the {\it adjoint relation} $A^+\in\widetilde\cC(\cH)$
is defined as the orthogonal companion of $A$ in $(\cH^2,\senkit)$, i.e.,
\begin{equation*}
A^+:=A^{\leftil\bot\rightil}=\bigl\{\hat f\in\cH^2:\leftik\hat f,\hat g\rightik=0\,\,\text{for all}\,\,
\hat g\in A\bigr\}.
\end{equation*}
A linear relation
$A$ in $\cH$ is said to be
\emph{symmetric} (\emph{selfadjoint}) if $A\subset A^*$ ($A=A^*$,
respectively).
We say that a closed symmetric relation $A\in\widetilde\cC(\cH)$ is of {\it defect} $m\in\dN_0\cup\{\infty\}$, if
the deficiency indices
\begin{equation*}
n_\pm(JA)=\dim\ker\bigl((JA)^*\mp i\bigr)
\end{equation*}
of the closed symmetric relation $JA$ in the Hilbert space
$(\cH,[J\cdot,\cdot])$ are both equal to $m$. Here $^*$ denotes the adjoint with respect
to the Hilbert scalar product $[J\cdot,\cdot]$. Note that a symmetric relation $A\in\widetilde\cC(\cH)$
is of defect $m$ if and only if there exists a selfadjoint
extension of $A$ in $\cH$ and each selfadjoint extension $A^\prime$
of $A$ in $\cH$ satisfies $\dim (A^\prime/A )=m$.

For symmetric operators in Hilbert spaces the concept of generalized boundary triples 
or generalized boundary value spaces
was introduced by V.A.~Derkach and M.M.~Malamud in \cite{DM95}, see also \cite[$\S 5.2$]{DHMS06}.
We use the same definition in the Krein space case.

\begin{definition}
Let $A$ be a closed symmetric relation in the Krein space $\cH$ and let $T$ be a linear relation in $\cH$ 
such that $\overline T=A^+$. A triple $\{\cG,\Gamma_0,\Gamma_1\}$ is said to be a {\em generalized boundary triple} 
for $A^+$, if $\cG$ is a Hilbert space 
and $\Gamma=(\Gamma_0,\Gamma_1)^\top:T\rightarrow\cG\times\cG$ is a linear mapping such that 
\begin{equation}\label{green}
\leftik\hat f,\hat g\rightik_{\cH^2}=\leftik\Gamma\hat f,\Gamma\hat g\rightik_{\cG^2}  
\end{equation}
holds for all $\hat f,\hat g\in T$, $\ran\Gamma_0=\cG$ and $A_0:=\ker\Gamma_0$ is a selfadjoint relation in $\cH$.
\end{definition}

Let $A\in\widetilde\cC(\cH)$ be a closed symmetric relation in $\cH$. 
Then a generalized boundary triple $\{\cG,\Gamma_0,\Gamma_1\}$ 
for $A^+$ exists if and only if $A$ admits a selfadjoint extension in $\cH$. In this case the defect of $A$
coincides with $\dim\cG$.
Assume now that $\{\cG,\Gamma_0,\Gamma_1\}$ is a generalized boundary triple for $A^+$.
Note that \eqref{green} can also be written in the form
\begin{equation}\label{greenexp}
 [f^\prime,g]-[f,g^\prime]=(\Gamma_1\hat f,\Gamma_0\hat g)_\cG-(\Gamma_0\hat f,\Gamma_1\hat g)_\cG,\quad \hat f=\{f,f^\prime\},\,
\hat g=\{g,g^\prime\}\in T.
\end{equation}
and that by \eqref{green} the operator $\Gamma:T\rightarrow\cG^2$, $T=\dom\Gamma$, is an isometry from the 
Krein space $(\cH^2,\senkit_{\cH^2})$ to the
Krein space $(\cG^2,\senkit_{\cG^2})$, i.e. $\Gamma^{-1}\subset\Gamma^{\senki+}$, where $\senki+$ denotes the adjoint
with respect to the Krein space inner products $\senkit_{\cH^2}$ in $\cH^2$ and $\senkit_{\cG^2}$ in $\cG^2$,
respectively. From
$\ran\Gamma_0=\cG$ and the selfadjointness of $A_0=\ker\Gamma_0$ one concludes that also the inclusion 
$\Gamma^{\senki+}\subset\Gamma^{-1}$ is true (cf. \cite[Lemma 5.5]{DHMS06}) and therefore
$\Gamma$ is a unitary operator from $(\cH^2,\senkit_{\cH^2})$ to $(\cG^2,\senkit_{\cG^2})$. This implies that $\Gamma$ is closed
and from \cite[Proposition 2.3]{DHMS06} we conclude  
$A=\ker\Gamma$ and that $\ran\Gamma$ is dense in $\cG^2$. Moreover, $\Gamma$ is surjective if and only if 
$\dom\Gamma=A^+$ holds.

Generalized boundary triples are a generalization of the well-known concept of (ordinary) 
boundary triples, see, e.g., \cite{BGP08,D95,D99,DM91,DM95,GG91}, and both notions coincide if 
the defect of the symmetric relation is finite. In short, a generalized boundary triple with 
a surjective $\Gamma$ is an ordinary boundary triple. The following definition from \cite{D99} reads
slightly different.

\begin{definition}\label{ordbt}
Let $A$ be a closed symmetric relation in the Krein space $\cH$. 
A triple $\{\cG,\Gamma_0,\Gamma_1\}$ is said to be an {\em ordinary boundary triple} 
for $A^+$, if $\cG$ is a Hilbert space 
and $\Gamma=(\Gamma_0,\Gamma_1)^\top:A^+\rightarrow\cG\times\cG$ is a surjective linear mapping such that 
\begin{equation}\label{greenabs}
\leftik\hat f,\hat g\rightik_{\cH^2}=\leftik\Gamma\hat f,\Gamma\hat g\rightik_{\cG^2}  
\end{equation}
holds for all $\hat f,\hat g\in A^+$.
\end{definition}

Let again $A\in\widetilde\cC(\cH)$ be symmetric and let $\{\cG,\Gamma_0,\Gamma_1\}$ be a 
generalized boundary triple for $A^+$, $T=\dom\Gamma$.
If the resolvent set $\rho(A_0)$ of the selfadjoint relation $A_0=\ker\Gamma_0$ is nonempty, then it is not
difficult to see that 
\begin{equation*}
A^+=A_0\,\widehat+\,\widehat\cN_{\lambda,A^+},\quad 
\widehat\cN_{\lambda,A^+}=\bigl\{\{f_\lambda,\lambda f_\lambda\}:f_\lambda\in\cN_{\lambda,A^+}=\ker(A^+-\lambda)\bigr\},
\end{equation*} 
holds for all $\lambda\in\rho(A_0)$. Here $\widehat +$ denotes the direct sum of subspaces. 
Since $\overline T=A^+$ and $A_0\subset T$ it follows that 
\begin{equation*}
\widehat\cN_{\lambda,T}=\bigl\{\{f_\lambda,\lambda f_\lambda\}:f_\lambda\in\cN_{\lambda,T}=\ker(T-\lambda)\bigr\}
\end{equation*}
is dense in $\widehat\cN_{\lambda,A^+}$ and $T$ can be decomposed as
\begin{equation}\label{dect}
T=A_0\,\widehat+\,\widehat\cN_{\lambda,T}=\ker\Gamma_0\,\widehat+\,\widehat\cN_{\lambda,T},\qquad\lambda\in\rho(A_0). 
\end{equation} 

Associated to a generalized boundary triple are the so-called $\gamma$-field and Weyl function.
For symmetric operators in Hilbert spaces the following definition can be found in \cite{DM95}.

\begin{definition}\label{gamweyldef}
Let $A$ be a closed symmetric relation in the Krein space $\cH$ 
and let $\{\cG,\Gamma_0,\Gamma_1\}$, $A_0=\ker\Gamma_0$, 
be a generalized boundary triple for $A^+$.
Assume $\rho(A_0)\not=\emptyset$ and denote the 
projection in $\cH\times\cH$ onto the first component by $\pi_1$. The {\em $\gamma$-field} $\gamma$ and
{\em Weyl function} $M$ corresponding to $\{\cG,\Gamma_0,\Gamma_1\}$ are defined by
\begin{equation*}
\gamma(\lambda)=\pi_1\bigl(\Gamma_0\!\upharpoonright\!\widehat\cN_{\lambda,T}\bigr)^{-1}\quad
\text{and}\quad M(\lambda)=\Gamma_1\bigl(\Gamma_0\!\upharpoonright\!\widehat\cN_{\lambda,T}\bigr)^{-1},
\quad\lambda\in\rho(A_0).
\end{equation*}
\end{definition}

In the following proposition we collect some properties of the $\gamma$-field and the Weyl function associated to
a generalized boundary triple. For $\gamma$-fields and Weyl functions of ordinary boundary triples
the statements in Proposition~\ref{prop1} are well known (see, e.g., \cite{D99})
and in our slightly more general situation the proofs are similar and in essence included in \cite[$\S$ 2.3]{BL07}.

\begin{proposition}\label{prop1}
Let $A\in\widetilde\cC(\cH)$ be symmetric, let 
$\{\cG,\Gamma_0,\Gamma_1\}$ be a generalized boundary triple for $A^+$
and assume $\rho(A_0)\not=\emptyset$, $A_0=\ker\Gamma_0$. Then the 
$\gamma$-field $\lambda\mapsto\gamma(\lambda)\in\cL(\cG,\cH)$ and Weyl function $\lambda\mapsto M(\lambda)\in\cL(\cG)$ of $\{\cG,\Gamma_0,\Gamma_1\}$ 
are holomorphic on $\rho(A_0)$ and the identities 
\begin{equation}\label{id1}
\gamma(\lambda)=\bigl(I+(\lambda-\mu)(A_0-\lambda)^{-1}\bigr)\gamma(\mu)
\end{equation}
and 
\begin{equation}\label{gambar}
\gamma(\bar\lambda)^+h=\Gamma_1\bigl\{(A_0-\lambda)^{-1}h,(I+\lambda(A_0-\lambda)^{-1})h\bigr\},\qquad h\in\cH,
\end{equation}
as well as 
\begin{equation}\label{id2}
M(\lambda)-M(\mu)^*=(\lambda-\bar\mu)\gamma(\mu)^+\gamma(\lambda)
\end{equation}
and 
\begin{displaymath}
M(\lambda)=\RE M(\lambda_0)+\gamma(\lambda_0)^+\bigl((\lambda-\RE\lambda_0)+(\lambda-\lambda_0)(\lambda-\bar\lambda_0)(A_0-\lambda)^{-1}\bigr)\gamma(\lambda_0)
\end{displaymath}
hold for all $\lambda,\mu\in\rho(A_0)$ and any fixed $\lambda_0\in\rho(A_0)$.
\end{proposition}

\section{Realization of operator functions as Weyl functions}\label{sec2}

Let $\cD\subset\dC^+$ be a simply connected open set, let $\cG$ be a Hilbert space and let 
$\tau$ be a piecewise meromorphic 
$\cL(\cG)$-valued function on $\cD\cup\cD^*$, $\cD^*=\{\lambda\in\dC:\bar\lambda\in\cD\}$,
which admits the representation
\begin{equation}\label{rep}
\tau(\lambda)=\RE\tau(\lambda_0)+\gamma^+\bigl((\lambda-\RE\lambda_0)+
(\lambda-\lambda_0)(\lambda-\bar\lambda_0)(A_0-\lambda)^{-1}\bigr)\gamma,
\end{equation}
with some selfadjoint relation $A_0$ in a Krein space $\cH$ and a mapping $\gamma\in\cL(\cG,\cH)$.
It is assumed that $\rho(A_0)$ is nonempty, that \eqref{rep} holds for a fixed 
$\lambda_0\in\cO\cup\cO^*$ and all
$\lambda\in\cO\cup\cO^*$, where $\cO$ is an open subset of $\rho(A_0)\cap\cD$, 
$\cO^*=\{\lambda\in\dC:\bar\lambda\in\cO\}$, and that the minimality condition
\begin{equation}\label{min}
\cH=\closp\bigl\{\bigl(I+(\lambda-\lambda_0)(A_0-\lambda)^{-1}\bigr)\gamma x:\lambda\in\cO\cup\cO^*,\,x\in\cG\bigr\}
\end{equation}
is satisfied. It is clear that $\tau$ is holomorphic on 
$\cO\cup\cO^*$ and that 
$\tau(\lambda)^*=\tau(\bar\lambda)$ holds for all $\lambda\in\cO\cup\cO^*$. 
The set of points of holomorphy of $\tau$ will be denoted by $\mathfrak h(\tau)$.

The following theorem is the main result of this section. The proof of Theorem~\ref{mainthm} will be given after 
some preparations at the end of in Section~\ref{sec2.2}. 

\begin{theorem}\label{mainthm}
Let $\tau:\cD\cup\cD^*\rightarrow\cL(\cG)$ be a piecewise meromorphic operator function which is represented
in the form \eqref{rep}-\eqref{min}. Then there exists a Krein space $\cK$, a closed symmetric
operator $S$ in $\cK$ and a generalized boundary triple 
$\{\cG,\Gamma_0,\Gamma_1\}$ for $S^+$ such
that the corresponding Weyl function coincides with $\tau$ on $\cO\cup\cO^*$. 
\end{theorem}

Since generalized boundary triples reduce to ordinary boundary triples if $\dim\cG$ is finite we obtain
the following corollary.

\begin{corollary}\label{mainthmcor}
Let $\tau:\cD\cup\cD^*\rightarrow\cL(\cG)$ be a piecewise meromorphic operator function which is represented
in the form \eqref{rep}-\eqref{min} and assume, in addition, that $\dim\cG$ is finite. 
Then there exists a Krein space $\cK$, a closed symmetric
operator $S$ in $\cK$ and an ordinary boundary triple 
$\{\cG,\Gamma_0,\Gamma_1\}$ for $S^+$ such
that the corresponding Weyl function coincides with $\tau$ on $\cO\cup\cO^*$. 
\end{corollary}

\begin{remark}\label{remrem}
Many important classes of $\cL(\cG)$-valued functions satisfy the above assumptions, cf. Section~\ref{sec2.3}. 
E.g., for 
Nevanlinna functions or generalized Nevanlinna functions  
one chooses $\cD=\dC^+$, $A_0$ becomes a selfadjoint relation in a Hilbert 
or Pontryagin space, respectively, and \eqref{rep} holds for all $\lambda\in\rho(A_0)$, cf. \cite{HSW98,KL77}.
So-called definitizable and locally
definitizable functions can be represented in the form \eqref{rep}-\eqref{min} 
with the help of definitizable and locally definitizable selfadjoint relations $A_0$ in Krein spaces, see 
\cite{J92,J00,J05}.
For operator functions piecewise holomorphic in $\cD\cup\cD^*$ and a given open subset $\cO$,
$\overline\cO\subset\cD$, a Krein space $\cH$ and a selfadjoint relation $A_0$ with 
$\overline\cO\cup\overline{\cO^*}\subset\rho(A_0)$ 
such that 
\eqref{rep}-\eqref{min} holds for all $\lambda\in\cO\cup\cO^*$ was constructed in \cite{A84,DLS87,J05}.
\end{remark}

Fix some $\mu_0\in\mathfrak h(\tau)$ and define the closed subspace $\widehat\cG$ of $\cG$ by 
\begin{equation}\label{stric}
\widehat\cG:=\bigcap_{\lambda\in\mathfrak h(\tau)}\ker\frac{\tau(\lambda)-\tau(\mu_0)^*}{\lambda-\bar\mu_0}.
\end{equation}
It is not difficult to see that $\widehat\cG$ does not depend on the choice of $\mu_0\in\mathfrak h(\tau)$ 
and that the set $\mathfrak h(\tau)$ in the intersection in \eqref{stric} 
can be replaced by the union of an open subset in 
$\cD$ and an open subset in $\cD^*$, e.g., $\cO\cup\cO^*$.

\begin{definition}\label{strictdef}
A piecewise meromorphic function $\tau:\cD\cup\cD^*\rightarrow\cL(\cG)$ is called {\em strict} if
the space $\widehat\cG$ in \eqref{stric} is trivial.
\end{definition}

\subsection{Realization of strict operator functions}

In this subsection we prove that every strict $\cL(\cG)$-valued operator function $\tau$ 
of the form \eqref{rep}-\eqref{min} can be realized as the Weyl function of a generalized boundary triple.
We start with a simple observation.

\begin{lemma}\label{lemstrict}
Let $\tau:\cD\cup\cD^*\rightarrow\cL(\cG)$ be a meromorphic function represented
in the form \eqref{rep}-\eqref{min} with some $\gamma\in\cL(\cG,\cH)$ and let $\widehat\cG$ be  
as in \eqref{stric}.
Then $\widehat\cG=\ker\gamma$ and, in particular, $\tau$ is strict if and only if $\gamma$ is injective.
\end{lemma}

\begin{proof}
For $x\in\ker\gamma$ we conclude from \eqref{rep}
$\tau(\lambda)x=\Real\tau(\lambda_0)x$ for all $\lambda\in\cO\cup\cO^*$ and therefore $x$ belongs to
\begin{equation}\label{cns}
\widehat\cG=\bigcap_{\lambda\in\mathfrak h(\tau)}\ker\,
\frac{\tau(\lambda)-\tau(\mu_0)^*}{\lambda-\bar\mu_0}.
\end{equation}
Conversely, if $x\in\widehat\cG$, then $x$ belongs also to the right hand side of 
\eqref{cns} with $\mu_0$ replaced by $\bar\lambda_0$. Making use of \eqref{rep} for $\lambda\in\cO\cup\cO^*$ 
we obtain
\begin{equation*}
\begin{split}
0&=\left(\frac{\tau(\lambda)-\tau(\lambda_0)}{\lambda-\lambda_0}\,x,y\right)
=\bigl(\gamma^{+}\bigl(I+(\lambda-\bar\lambda_0)(A_0-\lambda)^{-1}\bigr)
\gamma x,y\bigr)\\
&=\bigl[\gamma x,(I+(\bar\lambda-\lambda_0)
(A_0-\bar\lambda)^{-1})\gamma y\bigr]
\end{split}
\end{equation*}
for all $y\in\cG$ and all $\lambda\in\cO\cup\cO^*$.  The minimality condition \eqref{min}
implies $\gamma x=0$.
\end{proof}

The following theorem is a generalization of \cite[Theorem 3.3]{B08}, \cite[Proposition~3.1]{DHS01} and
\cite[$\S 3$]{DM95}.

\begin{theorem}\label{weylthm}
Let $\tau$ be a strict $\cL(\cG)$-valued function represented in the form \eqref{rep}-\eqref{min}. 
Then there exists a closed symmetric operator
$A$ in the Krein space $\cH$ and a generalized boundary triple 
$\{\cG,\Gamma_0,\Gamma_1\}$ for $A^+$ such that $\tau$ is the corresponding Weyl function on $\cO\cup\cO^*$.
Furthermore, $\{\cG,\Gamma_0,\Gamma_1\}$ is an
ordinary boundary triple if and only if $\ran\gamma$ is closed.
\end{theorem}

\begin{proof}
Let $\tau$ be represented by the selfadjoint relation $A_0$ in $\cH$
as in \eqref{rep}. For all
$\lambda\in\cO\cup\cO^*$ and the fixed
$\lambda_0\in\cO\cup\cO^*$ we define the mapping
\begin{equation}\label{gam}
\gamma(\lambda):=\bigl(I+(\lambda-\lambda_0)(A_0-\lambda)^{-1}\bigr)\gamma\;\in\cL(\cG,\cH).
\end{equation}
Then we have $\gamma(\lambda_0)=\gamma$,
$\gamma(\zeta)=(1+(\zeta-\eta)(A_0-\zeta)^{-1})\gamma(\eta)$
and 
\begin{equation}\label{qwer}
\tau(\zeta)-\tau(\eta)^*=(\zeta-\bar\eta)\gamma(\eta)^+\gamma(\zeta)
\end{equation}
for all
$\zeta,\eta\in\cO\cup\cO^*$.
For some $\xi\in\cO\cup\cO^*$ we define
the closed symmetric relation
\begin{equation}\label{tres}
A:=\bigl\{\{f_0,f_0^\prime\}\in A_0:[f_0^\prime-
\bar\xi f_0,\gamma(\xi) x]=0\;\;\text{for all}\;\;x\in\cG\bigr\}
\end{equation}
in $\cH$. Note that the definition of $A$ does not depend on the
choice of $\xi\in\cO\cup\cO^*$ and that
$\ran(A-\bar\lambda)=(\ran\gamma(\lambda))^{[\bot]}$ holds for all $\lambda\in\cO\cup\cO^*$.
Hence $\cN_{\lambda,A^+}=\overline{\ran\gamma(\lambda)}$
or, if $\ran\gamma(\lambda)$ is closed, then $\cN_{\lambda,A^+}=\ran\gamma(\lambda)$.
Since $\tau$ is assumed to be strict it follows from Lemma~\ref{lemstrict} that $\gamma$ 
is injective. Furthermore, the fact that the operator $I+(\lambda-\lambda_0)(A_0-\lambda)^{-1}$, 
$\lambda\in\cO\cup\cO^*$, is an isomorphism of
$\cN_{\lambda_0,A^+}$ onto $\cN_{\lambda,A^+}$ implies that
$\gamma(\lambda)$, regarded as a mapping from $\cG$ into
$\cN_{\lambda,A^+}$ is injective and has dense range.
Note also that the minimality condition~\eqref{min} together with \eqref{gam} 
implies that $A$ is an operator.

We fix a point $\mu\in\cO\cup\cO^*$. Then $A^+=A_0\,\widehat +\,\widehat\cN_{\mu,A^+}$
holds and the linear relation
\begin{equation*}
T:=A_0\,\widehat +\,\widehat\cN_{\mu,T},\qquad \widehat\cN_{\mu,T}=\bigl\{\{\gamma(\mu) x, \mu \gamma(\mu)x\}: x\in\cG\bigr\},
\end{equation*}
is dense in $A^+$. 
The elements
$\hat f\in T$
will be written in the form
\begin{equation*}
\hat f=\{ f_0, f_0^\prime\}+\{ \gamma(\mu) x,
\mu \gamma(\mu) x \},\quad \{f_0, f_0^\prime\}\in A_0,\, x\in\cG.
\end{equation*}
Let $\Gamma_0,\Gamma_1:T\rightarrow\cG$ be the linear mappings defined by
\begin{equation*}
\Gamma_0\hat f:=x\;\;\;\text{and}\;\;\;
\Gamma_1\hat f:=\gamma(\mu)^+(f_0^\prime-\bar{\mu}f_0)+\tau(\mu)x.
\end{equation*}
Then obviously $\ran\Gamma_0=\cG$ and $A_0=\ker\Gamma_0$ is selfadjoint. Moreover, for $\hat f\in T$ and
\begin{equation*}
\hat g=\{ g_0, g_0^\prime\}+\{\gamma(\mu) y,
\mu \gamma(\mu) y \}\in T,\quad \{ g_0, g_0^\prime\}\in A_0,\,y\in\cG,
\end{equation*}
we compute
\begin{displaymath}
\begin{split}
-i&\leftik\hat f,\hat g\rightik=[\gamma(\mu)x,g_0^\prime-\bar\mu g_0]-
[f_0^\prime-\bar\mu f_0,\gamma(\mu)y]-(\mu-\bar\mu)[\gamma(\mu)x,\gamma(\mu)y]\\
&=\bigl(x,\gamma(\mu)^+(g_0^\prime-\bar\mu g_0)\bigr)-\bigl(\gamma(\mu)^+(f_0^\prime-\bar\mu f_0),
y\bigr)
-\bigl((\tau(\mu)-\tau(\mu)^*)x,y\bigr)\\
&= -i\leftik\Gamma\hat f,\Gamma\hat g\rightik,
\end{split}
\end{displaymath}
where we have used $A_0=A_0^+$ and $\tau(\mu)-\tau(\mu)^*=(\mu-\bar\mu)\gamma(\mu)^+\gamma(\mu)$.
Therefore $\{\cG,\Gamma_0,\Gamma_1\}$ is a generalized boundary triple for $A^+$.

Let us check that the Weyl function corresponding to $\{\cG,\Gamma_0,\Gamma_1\}$
coincides with $\tau$ on $\cO\cup\cO^*$. Note first that by the definition of $\Gamma_0$ and $\Gamma_1$ it is clear
that $\tau(\mu)\Gamma_0\hat f_\mu=\Gamma_1\hat f_\mu$ holds for $\hat f_\mu=\{\gamma(\mu)x,\mu\gamma(\mu)x\}\in\widehat\cN_{\mu,T}$.
Now let $\eta\in\cO\cup\cO^*$ and $\hat f_\eta\in\hat\cN_{\eta,T}$. 
Since $T=A_0\,\widehat +\,\widehat\cN_{\mu,T}$
there exist $\{f_0,f_0^\prime\}\in A_0$
and $x\in\cG$ such that
\begin{equation}\label{etadeco}
\hat f_\eta=\{ f_\eta, \eta f_\eta\}
=\{ f_0, f_0^\prime\}+
\{\gamma(\mu)x , \mu\gamma(\mu) x \}.
\end{equation}
It follows from \eqref{qwer} and $\gamma(\eta)=(I+(\eta-\mu)(A_0-\eta)^{-1})\gamma(\mu)$ that 
\begin{equation*}
\begin{split}
\tau(\eta)&=\tau(\mu)^*+(\eta-\bar\mu)\gamma(\mu)^+\gamma(\eta)\\
&=\tau(\mu)+\gamma(\mu)^+\bigl((\bar\mu-\mu)\gamma(\mu)+
(\eta-\bar\mu)\gamma(\eta)\bigr)\\
&=\tau(\mu)+\gamma(\mu)^+(\eta-\mu)
\bigl(I+(\eta-\bar\mu)(A_0-\eta)^{-1}\bigr)\gamma(\mu).
\end{split}
\end{equation*}
Hence we have
\begin{equation}\label{taueta}
\tau(\eta)\Gamma_0\hat f_\eta=\tau(\mu)x
+\gamma(\mu)^+(\eta-\mu)
\bigl(I+(\eta-\bar\mu)(A_0-\eta)^{-1}\bigr)\gamma(\mu)x
\end{equation}
and from \eqref{etadeco} it follows that
\begin{equation*}
f_0^\prime-\eta f_0=(\eta-\mu)\gamma(\mu)x\quad\text{and}\quad
f_0^\prime-\bar\mu f_0=(\eta-\mu)\gamma(\mu)x+(\eta-\bar\mu)f_0
\end{equation*}
hold. The first identity yields $f_0=(\eta-\mu)(A_0-\eta)^{-1}\gamma(\mu)x$ and therefore \eqref{taueta} becomes
\begin{equation*}
\tau(\eta)\Gamma_0\hat f_\eta=\tau(\mu)x+\gamma(\mu)^+(f_0^\prime-\bar\mu f_0)
=\Gamma_1\hat f_\eta,
\end{equation*}
i.e., $\tau$ coincides with the Weyl function of $\{\cG,\Gamma_0,\Gamma_1\}$ on $\cO\cup\cO^*$.

It remains to show that the triple 
$\{\cG,\Gamma_0,\Gamma_1\}$ is an ordinary boundary triple for $A^+$ 
if and only if $\ran\gamma=\overline{\ran\gamma}$. Clearly, if $\{\cG,\Gamma_0,\Gamma_1\}$ is an ordinary boundary
triple, then the range of the $\gamma$-field is closed and hence $\ran\gamma=\ran\gamma(\lambda_0)$ is closed. 
Conversely, if $\ran\gamma$ is closed it is sufficient
to check that $(\Gamma_0,\Gamma_1)^\top$ is surjective, cf. Section~\ref{sec1}. Observe first 
that $\{0\}=\ker\gamma(\mu)=(\ran\gamma(\mu)^+)^\bot$ and that $\ran\gamma(\lambda)$ is closed for every $\lambda\in\cO\cup\cO^*$.
Hence $\ran\gamma(\mu)^+=\cG$ and for given elements $x,y\in\cG$ there exist $\{f_0,f_0^\prime\}\in A_0$ such that
$\gamma(\mu)^+(f_0^\prime-\bar\mu f_0)=y-\tau(\mu)x$. Now it easy to see that
$\hat f=\{f_0,f_0^\prime\}+\{\gamma(\mu)x,\mu\gamma(\mu)x\}$ satisfies $\Gamma_0\hat f=x$ 
and $\Gamma_1\hat f=y$.
\end{proof}

\begin{remark}\label{weylthmrem}
If $\tau$ is a strict $\cL(\cG)$-valued function which admits a representation as in \eqref{rep}-\eqref{min} and $\{\cG,\Gamma_0,\Gamma_1\}$
is a generalized boundary triple as in Theorem~\ref{weylthm} with $T=\dom\Gamma$, then the span of the
subspaces of $\cN_{\lambda,T}$ is dense in $\cH$, i.e., $\cH=\closp\{\cN_{\lambda,T}:\lambda\in\cO\cup\cO^*\}$, and the closed 
symmetric operator $A=\ker\Gamma$ has no eigenvalues.
\end{remark}

If $\tau$ is a matrix-valued function, that is, $\dim\cG<\infty$, then of course the range of
the mapping $\gamma\in\cL(\cG,\cH)$ in \eqref{rep} is closed. Hence Theorem~\ref{weylthm} implies 
the following corollary.

\begin{corollary}
Let $\tau$ be a strict $\cL(\cG)$-valued function represented in the form \eqref{rep}-\eqref{min}
and assume, in addition, that $\dim\cG$ is finite. 
Then there exists a closed symmetric operator
$A$ in the Krein space $\cH$ and an ordinary boundary triple 
$\{\cG,\Gamma_0,\Gamma_1\}$ for $A^+$ such that $\tau$ is the corresponding Weyl function on $\cO\cup\cO^*$.
\end{corollary}

\subsection{Realization of non-strict operator functions}\label{sec2.2}

Let again $\tau:\cD\cup\cD^*\rightarrow\cL(\cG)$ be a piecewise meromorphic operator function
which is represented in the form \eqref{rep}-\eqref{min}. We are now interested in the case 
where $\tau$ is not strict, i.e., the space $\widehat\cG$ in \eqref{stric} is not trivial. 
Roughly speaking the next lemma states that $\tau$ can always be written as a selfadjoint constant and a 
smaller strict operator function. For special classes of matrix-valued functions Lemma~\ref{reduce}
can be found in \cite{B08}.

\begin{lemma}\label{reduce}
Let $\tau$ be a piecewise meromorphic $\cL(\cG)$-valued function 
represented in the form \eqref{rep}-\eqref{min}, let $\widehat\cG$ be as in \eqref{stric}
and set $\cG^\prime:=\cG\ominus\widehat\cG$. Denote the corresponding orthogonal projections and canonical 
embeddings by $\widehat\pi$, $\pi^\prime$, $\widehat\iota$ and $\iota^\prime$, respectively, 
and fix some $\mu_0\in\mathfrak h(\tau)$. 
Then
\begin{equation}\label{taured}
\tau(\lambda)=\begin{pmatrix} \pi^\prime\tau(\lambda)\iota^\prime  & 0 \\ 0 & 0 \end{pmatrix}+
\begin{pmatrix} 0 & \pi^\prime\tau(\mu_0)\widehat\iota\\
\widehat\pi\tau(\mu_0)\iota^\prime & \widehat \pi\tau(\mu_0)\widehat\iota
\end{pmatrix}:
\begin{pmatrix}\cG^\prime\\ \widehat\cG\end{pmatrix}\rightarrow
\begin{pmatrix}\cG^\prime\\ \widehat\cG\end{pmatrix}
\end{equation}
for all $\lambda\in\mathfrak h(\tau)$ and the $\cL(\cG^\prime)$-valued
function $\lambda\mapsto \pi^\prime\tau(\lambda)\iota^\prime$ is strict.
\end{lemma}

\begin{proof}
It follows from the definition of $\widehat \cG$ in \eqref{stric} that
for $\widehat x\in\widehat\cG$ and all $\lambda\in\mathfrak h(\tau)$ 
the relation
$\tau(\lambda)\widehat\iota\widehat x
=\tau(\bar\mu_0)\widehat\iota\widehat x=\tau(\mu_0)\widehat\iota\widehat x$
holds. Therefore 
\begin{equation*}
\tau(\lambda)=\begin{pmatrix} \cdot & \pi^\prime\tau(\mu_0)\widehat\iota\\ 
\cdot & \widehat\pi\tau(\mu_0)\widehat\iota \end{pmatrix}:
\begin{pmatrix}\cG^\prime\\ \widehat\cG\end{pmatrix}\rightarrow
\begin{pmatrix}\cG^\prime\\ \widehat\cG\end{pmatrix},\quad\lambda\in\mathfrak h(\tau),
\end{equation*}
and the 
symmetry property $\tau(\bar\lambda)=\tau(\lambda)^*$
implies 
\begin{equation*}
\widehat\pi\tau(\lambda)\iota^\prime=(\pi^\prime\tau(\bar\lambda0\widehat\iota)^*=(\pi^\prime\tau(\bar\mu_0)\widehat\iota)^*=
\widehat\pi\tau(\mu_0)\iota^\prime
\end{equation*} 
which yields the representation \eqref{taured}.
Let us show that $\lambda\mapsto\pi^\prime\tau(\lambda)\iota^\prime$ is a strict function. Assume that 
$x^\prime\in\cG^\prime$ belongs to 
\begin{equation*}
\bigcap_{\lambda\in\mathfrak h(\tau)}\ker
\frac{\pi^\prime\tau(\lambda)\iota^\prime-\pi^\prime\tau(\bar\mu_0)\iota^\prime}{\lambda-\bar\mu_0}.
\end{equation*}
Then $\pi^\prime\tau(\lambda)\iota^\prime x^\prime=\pi^\prime\tau(\bar\mu_0)\iota^\prime x^\prime$ 
and also $\widehat\pi\tau(\lambda)\iota^\prime x^\prime=
\widehat\pi\tau(\bar\mu_0)\iota^\prime x^\prime$ by \eqref{taured} for all $\lambda\in\mathfrak h(\tau)$,
and this implies $\iota^\prime x^\prime\in\widehat\cG$. 
This is possible only for $x^\prime=0$, i.e., the function $\lambda\mapsto\pi^\prime\tau(\lambda)\iota^\prime$ is strict.
\end{proof}

Next we construct a nondensely defined closed symmetric operator $B$ in a Krein space
and an ordinary boundary triple for $B^+$ such that the corresponding Weyl function is a selfadjoint 
constant.

\begin{lemma}\label{consti}
Let $\widehat\cG$ be a Hilbert space, let $\Theta=\Theta^*\in\cL(\widehat\cG)$ and fix some $\vartheta\in\dC$.
Then $\widetilde\cH=(\widehat\cG^2,(J\cdot,\cdot))$, where 
$J=\bigl(\begin{smallmatrix} 0 & I \\ I & 0\end{smallmatrix}\bigr)$,
is a Krein space and there exists 
a closed symmetric operator $B$ in $\widetilde\cH$ and an ordinary boundary triple 
$\{\widehat\cG,\widehat\Gamma_0,\widehat\Gamma_1\}$, $B_0=\ker\widehat\Gamma_0$, 
for $B^+$ such that the corresponding Weyl function
is the selfadjoint constant $\Theta$ and $\sigma(B_0)=\{\vartheta,\bar\vartheta\}$.
\end{lemma}

\begin{proof}
We equip $\widehat\cG\times\widehat\cG$ with the indefinite inner product $[\cdot,\cdot]:=(J\cdot,\cdot)$, where
$J=\bigl(\begin{smallmatrix} 0 & I\\ I & 0\end{smallmatrix}\bigr)$ and $(\cdot,\cdot)$ is the
Hilbert scalar product on $\widehat\cG^2$. Then
\begin{equation*}
B_0:=\begin{pmatrix} \vartheta & I \\ 0 & \bar\vartheta \end{pmatrix}\in\cL(\widehat\cG^2)
\end{equation*}
is selfadjoint in the Krein space $\widetilde\cH=(\widehat\cG^2,[\cdot,\cdot])$ and
for every $\lambda\in\dC\backslash\{\vartheta,\bar\vartheta\}$ we have
\begin{equation*}
(B_0-\lambda)^{-1}=\begin{pmatrix} (\vartheta-\lambda)^{-1} & 
(\lambda-\vartheta)^{-1}(\bar\vartheta-\lambda)^{-1} 
\\ 0 & (\bar\vartheta-\lambda)^{-1}
\end{pmatrix}\in\cL(\widetilde\cH).
\end{equation*}
Let $\lambda_0\in\dC\backslash\{\vartheta,\bar\vartheta\}$,
$\widehat\gamma_{\lambda_0}:\widehat\cG\rightarrow\widetilde\cH$, $x\mapsto (x,0)^\top$, and define for
$\lambda\in\dC\backslash\{\vartheta,\bar\vartheta\}$
\begin{displaymath}
\widehat\gamma(\lambda): \widehat\cG\rightarrow\widetilde\cH,\qquad
x\mapsto\bigl(I+(\lambda-\lambda_0)(B_0-\lambda)^{-1}\bigr)\widehat\gamma_{\lambda_0} x
=\Bigl( \frac{\vartheta-\lambda_0}{\vartheta-\lambda}\,x, 0 \Bigr)^\top.
\end{displaymath}
Then obviously $\ran\widehat\gamma(\lambda)=\widehat\cG\times\{0\}$.
From
\begin{equation}\label{gammad}
\widehat\gamma(\eta)^+:\widetilde\cH\rightarrow\widehat\cG,\quad (x,y)^\top\mapsto 
\frac{\bar\vartheta-\bar\lambda_0}
{\bar\vartheta-\bar\eta}\, y,\quad
\eta\in\dC\backslash\{\vartheta,\bar\vartheta\},
\end{equation}
we obtain $\widehat\gamma(\eta)^+\widehat\gamma(\lambda)=0$ for all 
$\lambda,\eta\in\dC\backslash\{\vartheta,\bar\vartheta\}$. 
Consider
the closed symmetric operator
\begin{equation}\label{b0}
B:=B_0\upharpoonright\bigl(\widehat\cG\times\{0\}\bigr)
\end{equation}
in $\widetilde\cH$.
Then we have $\cN_{\lambda,B^+}=\widehat\cG\times\{0\}=\ran\widehat\gamma(\lambda)$ for
all $\lambda\in\dC\backslash\{\vartheta,\bar\vartheta\}$, the defect of $B$ coincides with 
$\dim\widehat\cG$ and
$\cN_{\lambda,B^+}[\bot]\cN_{\eta,B^+}$ holds for all $\lambda,\eta\in\dC
\backslash\{\vartheta,\bar\vartheta\}$.
For a fixed
$\mu\in\dC\backslash\{\vartheta,\bar\vartheta\}$
we write the elements $\hat g\in B^+=B_0\,\widehat +\,\widehat\cN_{\mu,B^+}$ in the form
\begin{equation*}
\hat g= \{ g_0, B_0 g_0\} +
\{ \widehat\gamma(\mu)x , \mu \widehat\gamma(\mu) x \},\qquad
g_0\in\widetilde\cH,\quad x\in\widehat\cG.
\end{equation*}
Then it follows as in the proof of Theorem~\ref{weylthm}
that $\{\widehat\cG,\widehat\Gamma_0,\widehat\Gamma_1\}$, where
\begin{equation}\label{tripelchen}
\widehat\Gamma_0\hat g:=x \quad\text{and}\quad 
\widehat\Gamma_1\hat g:=\widehat\gamma(\mu)^+(B_0-\bar\mu) g_0 +\Theta x,
\end{equation}
is a boundary triple for $B^+$ and the corresponding Weyl function is the selfadjoint 
constant $\Theta\in\cL(\widehat\cG)$. 
\end{proof}

\begin{remark}\label{negindex}
Note that the negative and the positive index of the Krein space $\widetilde\cH=(\widehat\cG^2,(J\cdot,\cdot))$ in 
Proposition~\ref{consti} coincides with $\dim\widehat\cG$, that is, 
\begin{equation*}
\dim\bigl(\ker(J-I)\bigr)=\dim\bigl(\ker(J+I)\bigr)=\dim\widehat\cG.
\end{equation*}
\end{remark}

\begin{proofspecial}
Let $\tau:\cD\cup\cD^*\rightarrow\cL(\cG)$ be a (in general non-strict) piecewise meromorphic
function which is represented in the form \eqref{rep}-\eqref{min} for a fixed $\lambda_0\in
\cO\cup\cO^*$ and all $\lambda\in\cO\cup\cO^*$. Let $\widehat\cG$ be as in \eqref{stric},
set $\cG^\prime=\cG\ominus\widehat\cG$ and decompose $\tau$ as in \eqref{taured}.

Then by Lemma~\ref{reduce} the piecewise meromorphic 
function 
\begin{equation*}
\tau_s:=\pi^\prime\tau\iota^\prime:\cD\cup\cD^*\rightarrow\cL(\cG^\prime)
\end{equation*} 
is strict. 
Setting $\gamma^\prime:=\gamma\iota^\prime\in\cL(\cG^\prime,\cH)$ it follows directly from \eqref{rep} that
\begin{equation*}
\tau_s(\lambda)=\RE\tau_s(\lambda_0)+\gamma^{\prime\,+}\bigl((\lambda-\RE\lambda_0)+
(\lambda-\lambda_0)(\lambda-\bar\lambda_0)(A_0-\lambda)^{-1}\bigr)\gamma^\prime
\end{equation*}
holds for a fixed $\lambda_0\in\cO\cup\cO^*$ and all $\lambda\in\cO\cup\cO^*$. 
Furthermore, \eqref{min} together with the fact $\widehat\cG=\ker\gamma$, cf. Lemma~\ref{lemstrict}, implies that 
the minimality condition 
\begin{equation*}
\cH=\closp\bigl\{\bigl(1+(\lambda-\lambda_0)(A_0-\lambda)^{-1}\bigr)\gamma^\prime x^\prime:\lambda\in\cO\cup\cO^*,
\,x^\prime\in\cG^\prime\bigr\}
\end{equation*}
is satisfied.  
Therefore we can apply Theorem~\ref{weylthm} to the function $\tau_s$, i.e., $\tau_s$ 
coincides on $\cO\cup\cO^*$ with the Weyl function corresponding to some closed symmetric operator $A\subset A_0$ 
in the Krein space $\cH$ and a generalized boundary triple $\{\cG^\prime,\Gamma_0^\prime,\Gamma_1^\prime\}$ 
for the adjoint $A^+$. Note that $A_0=\ker\Gamma_0^\prime$ and that $\dom\Gamma^\prime$, $\Gamma^\prime=(\Gamma_0^\prime,
\Gamma_1^\prime)^\top$, is dense in $A^+$.

According to 
Lemma~\ref{consti} there exists a Krein space $\widetilde\cH$, a closed symmetric operator $B$ in 
$\widetilde\cH$
and an ordinary boundary triple $\{\widehat\cG,\widehat\Gamma_0,\widehat\Gamma_1\}$ such that the corresponding 
Weyl function is the selfadjoint constant
\begin{equation*}
 \widehat \pi\tau(\mu_0)\widehat\iota \in\cL(\widehat\cG).
\end{equation*}
Moreover, the spectrum of the selfadjoint relation $B_0=\ker\widehat\Gamma_0$ consists of a pair 
of eigenvalues $\{\vartheta,\bar\vartheta\}$
and it is no restriction to assume that $\vartheta,\bar\vartheta\not\in\cO\cup\cO^*$ holds. 

In the following we consider the closed symmetric operator $S:=A\times B$ in the Krein space 
$\cK:=\cH\times\widetilde\cH$
and its adjoint $S^+=A^+\times B^+$. Note that $\dom\Gamma^\prime\times B^+$ is dense in  $S^+$.
The elements in $\dom\Gamma^\prime\times B^+$ will be denoted in the form $\{\hat f,\hat g\}$,
$\hat f\in \dom\Gamma^\prime$, $\hat g\in B^+$. We claim that $\{\cG,\Gamma_0,\Gamma_1\}$, where
\begin{equation*}
\Gamma_0\{\hat f,\hat g\}:=\begin{pmatrix}  \Gamma_0^\prime\hat f\\ \widehat\Gamma_0\hat g\end{pmatrix}
\quad\text{and}\quad
\Gamma_1\{\hat f,\hat g\}:=\begin{pmatrix} \Gamma_1^\prime\hat f+\pi^\prime\tau(\mu_0)\widehat\iota\,\widehat\Gamma_0\hat g
\\ \widehat\Gamma_1\hat g+\widehat\pi\tau(\mu_0)\iota^\prime\,\Gamma_0^\prime\hat f
\end{pmatrix},
\end{equation*}
$\{\hat f,\hat g\}\in \dom\Gamma^\prime\times B^+$, is a generalized boundary triple for $S^+$ such that 
the corresponding
Weyl function coincides with $\tau$ on $\cO\cup\cO^*$. In fact, since 
$\{\cG^\prime,\Gamma_0^\prime,\Gamma_1^\prime\}$ and
$\{\widehat\cG,\widehat\Gamma_0,\widehat\Gamma_1\}$ are generalized and ordinary boundary triples 
for $A^+$ and $B^+$, respectively,
it follows that for $\{\hat f,\hat g\},
\{\hat h,\hat k\}\in\dom\Gamma^\prime\times B^+$ 
\begin{displaymath}
\begin{split}
&\leftik\Gamma\{\hat f,\hat g\},\Gamma\{\hat h,\hat k\}\rightik_{(\cG^\prime\oplus\widehat\cG)^2}\\&=
i\left(\begin{pmatrix}  \Gamma_0^\prime\hat f\\ \widehat\Gamma_0\hat g\end{pmatrix},
\begin{pmatrix} \Gamma_1^\prime\hat h+\pi^\prime\tau(\mu_0)\widehat\iota\,\widehat\Gamma_0\hat k
\\ \widehat\Gamma_1\hat k+\widehat\pi\tau(\mu_0)\iota^\prime\,\Gamma_0^\prime\hat h
\end{pmatrix}
\right)
-i\left(\begin{pmatrix} \Gamma_1^\prime\hat f+\pi^\prime\tau(\mu_0)\widehat\iota\,\widehat\Gamma_0\hat g
\\ \widehat\Gamma_1\hat g+\widehat\pi\tau(\mu_0)\iota^\prime\,\Gamma_0^\prime\hat f
\end{pmatrix}, \begin{pmatrix}  \Gamma_0^\prime\hat h\\ \widehat\Gamma_0\hat k\end{pmatrix}\right)\\
&=\leftik\Gamma^\prime\hat f,\Gamma^\prime\hat h\rightik_{\cG^{\prime 2}}
+\leftik\widehat\Gamma\hat g,\widehat\Gamma\hat k\rightik_{\widehat\cG^{2}}
=\leftik\hat f,\hat h\rightik_{\cH^2}+\leftik\hat g,\hat k\rightik_{\widetilde\cH^2}
=\leftik\{\hat f,\hat g\},\{\hat h,\hat k\}\rightik_{(\cH\times\widetilde\cH)^2}
\end{split}
\end{displaymath}
holds. Here we also have used $(\pi^\prime\tau(\mu_0)\widehat\iota\,)^*=\widehat\pi\tau(\mu_0)\iota^\prime$.
Moreover, since $A_0=\ker\Gamma_0^\prime$ and $B_0=\ker\widehat\Gamma_0$ are selfadjoint in $\cH$ and $\widetilde\cH$,
respectively, it is clear that $\ker\Gamma_0=A_0\times B_0$ is a selfadjoint relation in 
$\cK=\cH\times\widetilde\cH$. As $\ran\Gamma_0^\prime=\cG^\prime$ and $\ran\widehat\Gamma_0=\widehat\cG$ we also have
that $\ran\Gamma_0$ coincides with $\cG=\cG^\prime\oplus\widehat\cG$. Hence $\{\cG,\Gamma_0,\Gamma_1\}$ 
is a generalized boundary triple for $S^+=A^+\times B^+$.
It remains to show that the corresponding Weyl function coincides with $\tau$. For this, note that
\begin{equation*}
\widehat\cN_{\lambda,\dom\Gamma}=\widehat\cN_{\lambda,\dom\Gamma^\prime\times B^+}
=\widehat\cN_{\lambda,\dom\Gamma^\prime}\times\widehat\cN_{\lambda,B^+},
\qquad\lambda\in\cO\cup\cO^*,
\end{equation*}
and let $\{\hat f_\lambda,\hat g_\lambda\}\in\dom\Gamma^\prime\times B^+$, where 
$\hat f_\lambda\in \widehat\cN_{\lambda,\dom\Gamma^\prime}$ and $\hat g_\lambda\in\widehat\cN_{\lambda,B^+}$. Since 
\begin{equation*}
\tau_s(\lambda)\Gamma_0^\prime\hat f_\lambda=
\Gamma_1^\prime\hat f_\lambda\quad\text{and}\quad
\widehat \pi\tau(\mu_0)\widehat\iota\,\widehat\Gamma_0\hat g_\lambda=
\widehat\Gamma_1\hat g_\lambda,\quad\lambda\in\cO\cup\cO^*, 
\end{equation*} 
we conclude
\begin{equation*}
\begin{split}
\tau(\lambda)\Gamma_0\{\hat f_\lambda,\hat g_\lambda\}&=\begin{pmatrix}\tau_s(\lambda) & \pi^\prime\tau(\mu_0)\widehat\iota\\
\widehat\pi\tau(\mu_0)\iota^\prime & \widehat \pi\tau(\mu_0)\widehat\iota\end{pmatrix}
\begin{pmatrix}  \Gamma_0^\prime\hat f_\lambda\\ \widehat\Gamma_0\hat g_\lambda\end{pmatrix}\\
&=
\begin{pmatrix} \Gamma_1^\prime\hat f_\lambda+ \pi^\prime\tau(\mu_0)\widehat\iota\,\widehat\Gamma_0\hat g_\lambda\\
\widehat\pi\tau(\mu_0)\iota^\prime\,\Gamma_0^\prime\hat f_\lambda+\widehat\Gamma_1\hat g_\lambda
\end{pmatrix}=\Gamma_1\{\hat f_\lambda,\hat g_\lambda\}
\end{split}
\end{equation*}
for all $\lambda\in\cO\cup\cO^*$, that is, $\tau$ coincides with the Weyl function corresponding to $\{\cG,\Gamma_0,\Gamma_1\}$
on $\cO\cup\cO^*$.
\end{proofspecial}

\begin{remark}\label{nonminrem}
Let $\tau$ be as in \eqref{rep}-\eqref{min} and let $\cK=\cH\times\widetilde\cH$, $S=A\times B$ and 
$\{\cG,\Gamma_0,\Gamma_1\}$
be as in the proof of Theorem~\ref{mainthm}. If $\tau$ is non-strict, then $\widehat\cG\not=\{0\}$ and 
in contrast to Theorem~\ref{weylthm} and 
Remark~\ref{weylthmrem} here the defect subspaces $\cN_{\lambda,\dom\Gamma}$, $\lambda\in\cO\cup\cO^*$, 
are not dense in $\cK$. Indeed, it follows from the construction in the proof of
Lemma~\ref{consti} that
\begin{equation*}
\closp\bigl\{\cN_{\lambda,B^+}:\lambda\in\cO\cup\cO^*\bigr\}=\widehat\cG\times\{0\}\not=
\widetilde\cH=\widehat\cG\times\widehat\cG
\end{equation*}
holds. Therefore
\begin{equation*}
\closp\bigl\{\cN_{\lambda,\dom\Gamma}:\lambda\in\cO\times\cO^*\bigr\}=\cH\times\widehat\cG\times\{0\}\not=
\cK.
\end{equation*}
This implies that the analytic properties of $\tau$ are in general not completely reflected by the spectral 
properties of the selfadjoint operator or relation $S_0=\ker\Gamma_0$ in $\cK$, but this disadvantage arises
only at the points $\vartheta, \bar\vartheta$ which can be chosen arbitrary, e.g. in $\dC\backslash (\cD\cup\cD^*)$. In Section~\ref{pde} 
we shall see that the non-minimality does not affect solvability properties of a certain class of elliptic boundary
value problems investigated here. Note also, that $\vartheta$ is the only eigenvalue of the symmetric operator
$S=A\times B$, since $\sigma_p(A)=\emptyset$ by Remark~\ref{weylthmrem} and $\sigma_p(B)=\{\vartheta\}$; cf. \eqref{b0}.
\end{remark}

\subsection{Some special classes of operator functions}\label{sec2.3}

Many classes of $\dR$-symmetric operator functions satisfy the general assumptions in the beginning of
Section~\ref{sec2}, cf. Remark~\ref{remrem}. In this subsection we briefly recall some necessary 
definitions and we formulate some corollaries of Theorem~\ref{mainthm}. 

The first corollary concerns the case of a locally holomorphic operator function. We refer to \cite{A84,DLS87,J05}
for the existence of the representation \eqref{rep}-\eqref{min}.

\begin{corollary}
Let $\tau:\cD\cup\cD^*\rightarrow\cL(\cG)$ be a piecewise holomorphic function which satisfies 
$\tau(\bar\lambda)=\tau(\lambda)^*$, $\lambda\in\cD\cup\cD^*$, and let $\cO$ be a simply
connected open set with $\overline\cO\subset\cD$. 
Then there exists a Krein space $\cK$, a closed symmetric
operator $S$ in $\cK$ and a generalized boundary triple 
$\{\cG,\Gamma_0,\Gamma_1\}$ for $S^+$ such
that the corresponding Weyl function coincides with $\tau$ on $\cO\cup\cO^*$. If, in addition, $\dim\cG<\infty$
holds, then $\{\cG,\Gamma_0,\Gamma_1\}$ is an ordinary boundary triple.
\end{corollary}

The classes of generalized Nevanlinna functions were introduced and studied by M.G.~Krein and H.~Langer, 
see, e.g., \cite{KL72,KL77,KL81}.
Recall that an $\cL(\cG)$-valued function $\tau$ belongs to the {\it generalized Nevanlinna class} 
$N_\kappa(\cL(\cG))$,
$\kappa\in\dN_0$, 
if $\tau$ is piecewise meromorphic in $\dC\backslash\dR$ and $\dR$-symmetric, i.e., 
$\tau(\bar\lambda)=\tau(\lambda)^*$ for all $\lambda$ belonging to the set of points of holomorphy $\mathfrak h(\tau)$ of $\tau$, 
and the kernel
\begin{equation*}
K_\tau(\lambda,\mu):=\frac{\tau(\lambda)-\tau(\mu)^*}{\lambda-\bar\mu},\qquad \lambda,\mu\in\dC^+\cap
\mathfrak h(\tau),
\end{equation*}
has $\kappa$ negative squares, that is, for all $n\in\dN$, $\lambda_1,\dots,\lambda_n\in\dC^+\cap
\mathfrak h(\tau)$ and all $x_1,\dots,x_n\in\cG$ the selfadjoint matrix
\begin{equation*}
\bigl((K_\tau(\lambda_i,\lambda_j)x_i,x_j)\bigl)_{i,j=1}^n
\end{equation*}
has at most $\kappa$ negative eigenvalues, and $\kappa$ is minimal with this property. The functions
in the class $N_0(\cL(\cG))$ are called {\it Nevanlinna functions}. A function $\tau\in N_0(\cL(\cG))$ is holomorphic
on $\dC\backslash\dR$ and $\Imag\tau(\lambda)$ is nonnegative for all $\lambda\in\dC^+$.  
It is well-known that Nevanlinna functions can equivalently be characterized by integral representations. 
More precisely, 
$\tau$ is a $\cL(\cG)$-valued Nevanlinna function if and only if there exist selfadjoint operators $\alpha,\beta\in\cL(\cG)$,
$\beta\geq 0$, and a nondecreasing selfadjoint operator function $t\mapsto\Sigma(t)\in \cL(\cG)$ on $\dR$
such that 
$\int_\dR\tfrac{1}{1+t^2}d\Sigma(t)\in\cL(\cG)$ and
\begin{equation}\label{intrep}
\tau(\lambda)=\alpha+\lambda
\beta+\int_{-\infty}^\infty\Big(\frac1{t-\lambda}-\frac{t}{1+t^2}\Big)
\,d\Sigma(t)
\end{equation}
holds for all $\lambda\in\mathfrak h(\tau)$. It is worth to note that a Nevanlinna function $\tau$ is 
strict if and only if $\Imag\tau(\lambda)$ is uniformly positive for 
some (and hence for all) $\lambda\in\dC^+$.

It was shown in \cite{HSW98,KL77} that every function $\tau\in N_\kappa(\cL(\cG))$ can be represented 
in the form \eqref{rep}-\eqref{min} with $\cD=\dC^+$, $\cO=\mathfrak h(\tau)\cap\dC^+$ and $\cH$ is
a Pontryagin space with negative index $\kappa$. For generalized Nevanlinna functions our main result
reads as follows, cf. Remark~\ref{negindex} and \cite[Theorem 3.2]{PAMS} for the special case of $\cL(\dC^n)$-valued
Nevanlinna functions.

\begin{corollary}\label{nkappacor}
Let $\tau\in N_\kappa(\cL(\cG))$, $\kappa\in\dN_0$, and let $\widehat\cG$ be as in \eqref{stric}. Then 
there exists a Krein space $\cK$ with negative index $\kappa+\dim\widehat\cG$, a closed symmetric
operator $S$ in $\cK$ and a generalized boundary triple 
$\{\cG,\Gamma_0,\Gamma_1\}$ for $S^+$ such
that the corresponding Weyl function coincides with $\tau$ on $\mathfrak h(\tau)$.
If, in addition, $\dim\cG<\infty$, then $\cK$ is a Pontryagin space 
with negative index $\kappa+\dim\widehat\cG$ and $\{\cG,\Gamma_0,\Gamma_1\}$ is an ordinary boundary triple. 
\end{corollary}

Next we briefly recall the definitions of definitizable and locally definitizable operator functions
introduced by P. Jonas in \cite{J92,J00,J05}. An $\dR$-symmetric 
piecewise meromorphic $\cL(\cG)$-valued function $\tau$ in $\dC\backslash\dR$ is called {\it definitizable}
if there exists an $\dR$-symmetric scalar rational function $r$ such that $r\tau$ is the sum of a 
Nevanlinna function $G\in N_0(\cL(\cG))$ and an $\cL(\cG)$-valued rational function $P$
with the poles of $P$ belonging to $\mathfrak h(\tau)$,
\begin{equation*}
r(\lambda)\tau(\lambda)=G(\lambda)+P(\lambda),\qquad \lambda\in\mathfrak h(r\tau).
\end{equation*}
The classes $N_\kappa(\cL(\cG))$, $\kappa\in\dN_0$, are contained in the set of definitizable functions, 
see \cite{J92,J00}. 
Let $\Omega$ be a domain in $\overline\dC$ which is symmetric with respect to $\dR$, such that 
$\Omega\cap\overline\dR\not=\emptyset$ and $\Omega\cap\dC^+$ and $\Omega\cap\dC^-$ are simply connected.
A $\cL(\cG)$-valued function $\tau$ is said to be {\it definitizable in} $\Omega$ if for every domain 
$\Omega^\prime$ with the same properties as $\Omega$, $\overline{\Omega^\prime}\subset\Omega$, the 
restriction of $\tau$ to $\Omega^\prime$ can be written as the sum of a definitizable function $\tau_d$ and
an $\dR$-symmetric $\cL(\cG)$-valued function $\tau_{h}$ holomorphic in $\Omega^\prime$, 
$\tau(\lambda)=\tau_d(\lambda)+\tau_{h}(\lambda)$ for all $\lambda\in\mathfrak h(\tau)\cap\Omega^\prime$.

Operator representations of the form \eqref{rep}-\eqref{min} 
for definitizable and locally definitizable functions can be found in \cite{J00,J05}.
If $\tau$ is definitizable in $\Omega$ and $\Omega^\prime$ is a domain as $\Omega$, 
$\overline{\Omega^\prime}\subset\Omega$, 
one can choose $\cD=\Omega\cap\dC^+$ and $\cO=\Omega^\prime\cap\mathfrak h(\tau)\cap\dC^+$.
This yields the following corollary.

\begin{corollary}
Let $\tau$ be a $\cL(\cG)$-valued function definitizable in $\Omega$ and let $\Omega^\prime$ be a domain 
with the same properties as $\Omega$, $\overline{\Omega^\prime}\subset\Omega$. 
Then there exists a Krein space $\cK$, a closed symmetric
operator $S$ in $\cK$ and a generalized boundary triple 
$\{\cG,\Gamma_0,\Gamma_1\}$ for $S^+$ such
that the corresponding Weyl function coincides with $\tau$ on 
$\Omega^\prime\cap\mathfrak h(\tau)\cap\dC\backslash\dR$. If, in addition, $\dim\cG<\infty$
holds, then $\{\cG,\Gamma_0,\Gamma_1\}$ is an ordinary boundary triple.
\end{corollary}

\section{Elliptic PDEs with $\lambda$-dependent boundary conditions}\label{pde}

Let $\Omega$ be a smooth bounded domain in $\dR^n$, $n>1$, with $C^\infty$-boundary $\partial\Omega$ and consider the 
second order differential expression
\begin{equation}\label{diffexp}
\ell =-\sum_{j,k=1}^n\partial_j\,a_{jk}\,\partial_k \,+a
\end{equation}
on $\Omega$ with coefficients $a_{jk},a\in C^\infty(\overline\Omega)$ such that
$a_{jk}=\overline{a_{kj}}$ for all $j,k=1,\dots,n$ and $a$ is real-valued. In addition, it is assumed that the ellipticity condition
\begin{equation*}
\sum_{j,k=1}^n a_{jk}(x)\xi_j\xi_k\geq C\sum_{k=1}^n \xi_k^2,\qquad \xi=(\xi_1,\dots,\xi_n)^\top\in\dR^n,\,\,x\in\overline\Omega,
\end{equation*}
holds for some constant $C>0$. In this section we investigate the following $\lambda$-dependent elliptic boundary value problem:
For a given function
$g\in L^2(\Omega)$ and $\lambda\in\mathfrak h(\tau)$ find $f\in L^2(\Omega)$ such that
\begin{equation}\label{rwp}
(\ell -\lambda) f=g\qquad\text{and}\qquad \tau(\lambda) f\vert_{\partial\Omega}=\frac{\partial f_D}{\partial\nu_\ell}\bigl\vert_{\partial\Omega}
\end{equation}
holds. Here $\tau$ is assumed to be a piecewise meromorphic $\cL(L^2(\partial\Omega))$-valued function and 
$f_D$ denotes the component of $f$ in the domain of the Dirichlet operator. The precise formulation of the problem will
be given in Section~\ref{bvpsec}.

\subsection{Preliminaries and ordinary boundary triples for elliptic PDEs}\label{pde1}

The Sobolev space of $k$th order on $\Omega$ is denoted by $H^k(\Omega)$ and the closure 
of $C_0^\infty(\Omega)$ in $H^k(\Omega)$ is denoted by $H_0^k(\Omega)$. Sobolev spaces on the boundary are denoted 
by $H^s(\partial\Omega)$, $s\in\dR$. Let $(\cdot,\cdot)_{-1/2\times 1/2}$ and $(\cdot,\cdot)_{-3/2\times 3/2}$ be the
extensions of the $L^2(\partial\Omega)$ inner product to $H^{-1/2}(\partial\Omega)\times H^{1/2}(\partial\Omega)$ and
$H^{-3/2}(\partial\Omega)\times H^{3/2}(\partial\Omega)$, respectively, and let $\iota_\pm : H^{\pm 1/2}(\partial\Omega)
\rightarrow L^2(\partial\Omega)$ be isomorphisms such that $(x,y)_{-1/2\times 1/2}=(\iota_- x,\iota_+ y)$ holds for all
$x\in H^{-1/2}(\partial\Omega)$ and $y\in H^{1/2}(\partial\Omega)$.

Recall that the {\it Dirichlet operator} 
\begin{equation*}
T_D f_D=\ell f_D,\qquad \dom T_D=H^2(\Omega)\cap H^1_0(\Omega),
\end{equation*}
associated to the elliptic differential expression $\ell$ in \eqref{diffexp} is selfadjoint in $L^2(\Omega)$ and the 
resolvent of $T_D$ is compact, cf. \cite[VI. Theorem~1.4]{EE87} and \cite{LU,LM72,W}.
Furthermore, the {\it minimal operator} 
\begin{equation*}
Tf=\ell f,\qquad \dom T=H^2_0(\Omega),
\end{equation*} 
is a densely defined 
closed symmetric operator in $L^2(\Omega)$ and the adjoint operator $T^*f=\ell f$ is defined on the maximal domain 
\begin{equation*}
\dom T^*=\cD_\max=\bigl\{f\in L^2(\Omega):\ell f\in L^2(\Omega)\bigr\}.
\end{equation*} 
Let us fix some
$\eta\in\dR\cap\rho(T_D)$. Then for each function $f\in\cD_\max$
there is a unique decomposition $f=f_D+f_\eta$, where $f_D\in\dom T_D$ and $f_\eta\in\cN_{\eta,T^*}=\ker(T^*-\eta)$.
In fact, as $T_D-\eta$ is surjective for a given $f\in\cD_\max$ there exists $f_D\in\dom T_D$ such that $(T^*-\eta)f=(T_D-\eta)f_D$
holds. It follows that $f_\eta:=f-f_D\in\cN_{\eta,T^*}$ and hence $f=f_D+f_\eta$ is the desired decomposition. 
The uniqueness follows from $\ker(T_D-\eta)=\{0\}$.

Let $\mathfrak n=(\mathfrak n_1,\dots,\mathfrak n_n)^\top$ be the unit outward normal of $\Omega$.
It is well-known that the map 
\begin{displaymath}
C^\infty(\overline\Omega)\ni f\mapsto \left\{f\vert_{\partial\Omega},\frac{\partial f}{\partial\nu_\ell}\Bigl|_{\partial\Omega}\right\},
\quad\text{where}\quad\frac{\partial f}{\partial\nu_\ell}:=\sum_{j,k=1}^n a_{jk} \mathfrak n_j \partial_k f,
\end{displaymath}
can be extended to a linear operator from $\cD_\max$ into $H^{-1/2}(\partial\Omega)\times H^{-3/2}(\partial\Omega)$ and that
for $f\in\cD_\max$ and $g\in H^2(\Omega)$ Green's identity
\begin{equation}\label{green2}
(T^*f,g)-(f,T^*g)=\left(f\vert_{\partial\Omega},\frac{\!\!\partial g}{\partial\nu_\ell}\bigl|_{\partial\Omega}\right)_{-\frac{1}{2}\times\frac{1}{2}}-
\left(\frac{\partial f}{\partial\nu_\ell}\bigl|_{\partial\Omega},g\vert_{\partial\Omega}\right)_{\!\!-\frac{3}{2}\times\frac{3}{2}}
\end{equation}
holds, see \cite{G68,LM72,W}.

The $\lambda$-dependent boundary condition in \eqref{rwp} will be rewritten with the help of an ordinary boundary triple for
the maximal realization of $\ell$ in $L^2(\Omega)$. The ordinary boundary triple in the next proposition can also be found 
in \cite{BGW09,GM08,GBook08,G08}.
For the convenience of 
the reader we include a short proof based on the general observations in \cite{G68,G71}.

\begin{proposition}\label{elbt}
The triple $\{L^2(\partial\Omega),\Upsilon_0,\Upsilon_1\}$, where 
\begin{equation*}
\Upsilon_0 \hat f:=\iota_- f_\eta\vert_{\partial\Omega}\qquad\text{and}\qquad 
\Upsilon_1 \hat f:=- \iota_+ \frac{\partial f_D}{\partial\nu_\ell}\bigl\vert_{\partial\Omega},
\end{equation*}
$\hat f=\{f,T^* f\}$, $f=f_D+f_\eta\in\cD_\max$,
is an ordinary boundary triple for the maximal operator $T^* f=\ell f$, $\dom T^*=\cD_\max$, such that $T_D=\ker\Upsilon_0$.
The corresponding $\gamma$-field and Weyl function
are given by
\begin{equation*}
\gamma(\lambda)y=(I+(\lambda-\eta)(T_D-\lambda)^{-1})f_\eta(y),\qquad\lambda\in\rho(T_D),
\end{equation*}
and 
\begin{equation*} 
M(\lambda)y=(\eta-\lambda)\iota_+
\frac{\partial (T_D-\lambda)^{-1}f_\eta(y)}{\partial \nu_\ell}\bigl|_{\partial\Omega},\qquad\lambda\in\rho(T_D),
\end{equation*}
respectively, where $f_\eta(y)$ is the unique function in $\ker(T^*-\eta)$ satisfying $\iota_- f_\eta(y)\vert_{\partial\Omega}=y$.
\end{proposition}

\begin{proof}
Let $f,g\in\cD_\max$ be decomposed in the form $f=f_D+f_\eta$ and $g=g_D+g_\eta$. As $T_D$ is selfadjoint and $\eta\in\dR$ we find
\begin{displaymath}
(T^*f,g)-(f,T^*g)=(T_Df_D,g_\eta)-(f_D,T^*g_\eta)+(T^*f_\eta,g_D)-(f_\eta,T_Dg_D)
\end{displaymath}
and then $f_D\vert_{\partial\Omega}=g_D\vert_{\partial\Omega}=0$ together with Green's identity \eqref{green2} implies
\begin{equation*}
\begin{split}
(T^*f,g)-(f,T^*g)&=-
\left(\frac{\partial f_D}{\partial\nu_\ell}\Bigl|_{\partial\Omega},g_\eta\vert_{\partial\Omega}\right)_{\frac{1}{2}\times -\frac{1}{2}}+
\left(f_\eta\vert_{\partial\Omega},\frac{\partial g_D}{\partial\nu_\ell}\Bigl|_{\partial\Omega}\right)_{-\frac{1}{2}\times \frac{1}{2}}\\
&=(\Upsilon_1 \hat f,\Upsilon_0\hat g)-(\Upsilon_0\hat f,\Upsilon_1\hat g).
\end{split}
\end{equation*}
Hence \eqref{greenabs} in Definition~\ref{ordbt} holds, cf. \eqref{greenexp}.
Furthermore, by the classical trace theorem the map 
$H^2(\Omega)\cap H^1_0(\Omega)\ni f_D\mapsto \frac{\partial f_D}{\partial\nu_\ell}\vert_{\partial\Omega}\in H^{1/2}(\partial\Omega)$
is onto and the same holds for the map $\ker(T^*-\eta)\ni f_\eta\mapsto f_\eta\vert_{\partial\Omega}\in H^{-1/2}(\partial\Omega)$, 
which is an isomorphism according to \cite[Theorem 2.1]{G71}. Hence $(\Upsilon_0,\Upsilon_1)^\top$ maps $T^*$ onto $L^2(\partial\Omega)\times L^2(\partial\Omega)$
and therefore $\{L^2(\partial\Omega),\Upsilon_0,\Upsilon_1\}$ is an ordinary boundary triple for $T^*$ with $T_D=\ker\Upsilon_0$.

It remains to show that the corresponding $\gamma$-field and Weyl function have the asserted form. For this let 
$y\in L^2(\partial\Omega)$, choose the unique function $f_\eta(y)$ in $\ker(T^*-\eta)$ such that 
$y=\iota_- f_\eta(y)\vert_{\partial\Omega}$ holds and set
\begin{equation}\label{flambda}
 f_\lambda:=(\lambda-\eta)(T_D-\lambda)^{-1}f_\eta(y)+f_\eta(y)
\end{equation}
for $\lambda\in\rho(T_D)$.
It is easy to see that $(T^*-\lambda)f_\lambda=0$ holds and since $(T_D-\lambda)^{-1}f_\eta(y)\in\dom T_D$ and 
$f_\eta(y)\in \ker(T^*-\eta)$ we obtain
\begin{equation*}
\Gamma_0\hat f_\lambda=\Gamma_0\{f_\lambda,\lambda f_\lambda\}=\iota_-f_\eta(y)\vert_{\partial\Omega}=y,
\end{equation*}
i.e. $\gamma(\lambda)y=f_\lambda=(I+(\lambda-\eta)(T_D-\lambda)^{-1})f_\eta(y)$. Finally, by the definition of the Weyl function 
and \eqref{flambda} we have
\begin{equation*}
M(\lambda) y=\Gamma_1 \hat f_\lambda=(\eta-\lambda)\iota_+ 
\frac{\partial (T_D-\lambda)^{-1}f_\eta(y)}{\partial \nu_\ell}\bigl|_{\partial\Omega}.
\end{equation*}
\end{proof}

\subsection{Elliptic boundary value problems with eigenvalue depending boundary conditions}\label{bvpsec}

Let $\cD\subset\dC^+$ be a simply connected open set and let $\tau$ be a piecewise meromorphic $\cL(L^2(\partial\Omega))$-valued 
function on $\cD\cup\cD^*$ which admits a representation of the form \eqref{rep}-\eqref{min} via the resolvent of some selfadjoint
relation on an open subset $\cO\cup\cO^*$ of $\cD\cup\cD^*$. Note that $\tau$ is holomorphic on $\cO\cup\cO^*$.
We study the following $\lambda$-dependent elliptic boundary value problem: For a given function
$g\in L^2(\Omega)$ and $\lambda\in\cO\cup\cO^*$ find $f\in\cD_\max$ such that
\begin{equation}\label{rwp1}
(\ell -\lambda) f=g\qquad\text{and}\qquad \tau(\lambda) \iota_- f\vert_{\partial\Omega}=\iota_+\frac{\partial f_D}{\partial\nu_\ell}\bigl\vert_{\partial\Omega}
\end{equation}
holds. According to Theorem~\ref{mainthm} there exists a Krein space $\cK$, a closed symmetric operator $S$ in  $\cK$ and a 
generalized boundary triple $\{L^2(\partial\Omega),\Gamma_0,\Gamma_1\}$ for $S^+=\overline{\dom\Gamma}$ such that the corresponding Weyl function
coincides with $\tau$ on $\cO\cup\cO^*$. In particular, the set $\cO\cup\cO^*$ is a subset of the resolvent 
set of the selfadjoint relation $S_0=\ker\Gamma_0$ in $\cK$.
With the help of the operator $S$, the generalized boundary triple $\{L^2(\partial\Omega),\Gamma_0,\Gamma_1\}$ 
and the ordinary boundary
triple $\{L^2(\partial\Omega),\Upsilon_0,\Upsilon_1\}$ for the elliptic operator from Proposition~\ref{elbt} we construct a linearization of the boundary value problem
\eqref{rwp1} in the next theorem.

\begin{theorem}\label{rwpthm}
Let $\{L^2(\partial\Omega),\Upsilon_0,\Upsilon_1\}$ be the ordinary boundary triple for the maximal differential operator 
$T^*$ associated to $\ell$ from Proposition~\ref{elbt}
with corresponding $\gamma$-field $\gamma$ and Weyl function $M$, and assume that 
$(M(\mu)+\tau(\mu))^{-1}\in\cL(L^2(\partial\Omega))$ holds for some $\mu\in\cO$.

Then the operator 
\begin{displaymath}
\begin{split}
\widetilde A\begin{pmatrix} f\\ k\end{pmatrix}&=\begin{pmatrix} \ell f \\ k^\prime\end{pmatrix},\\
\dom\widetilde A&=\left\{\begin{pmatrix} f \\ k \end{pmatrix}\in\cD_\max\times \cK:
\begin{pmatrix} \Upsilon_0\hat f - \Gamma_0\hat k  = 0\\ \Upsilon_1\hat f + \Gamma_1\hat k = 0\end{pmatrix}\,\,\text{for}\,\,\,\,
\begin{matrix}\hat f=\{f,T^*f\},\\ \hat k=\{k,k^\prime\}\in\dom\Gamma\end{matrix}\right\},
\end{split}
\end{displaymath}
is a selfadjoint extension of the minimal differential operator $T$ in the Krein space $L^2(\Omega)\times\cK$, the set 
\begin{equation*}
\cU:=\bigl\{\lambda\in\cO\cup\cO^*:(M(\lambda)+\tau(\lambda))^{-1}\in\cL(L^2(\partial\Omega))\bigr\}
\end{equation*}
is a subset of $\rho(\widetilde A)\cap\rho(T_D)\cap\mathfrak h(\tau)$ and for every $\lambda\in\cU$ 
the unique solution of 
the boundary value problem \eqref{rwp1} is given by
\begin{equation}\label{sol}
f=P_{L^2}\bigl(\widetilde A-\lambda)^{-1}\!\upharpoonright_{L^2} g=(T_D-\lambda)^{-1}g-\gamma(\lambda)\bigl(M(\lambda)+\tau(\lambda)\bigr)^{-1}
\gamma(\bar\lambda)^*g.
\end{equation}
\end{theorem}

\begin{proof}
The proof of Theorem~\ref{rwpthm} is divided into two parts. In the first part it will be shown that $\widetilde A$ is a selfadjoint
operator in the Krein space $L^2(\Omega)\times\cK$ and in the second part it is verified that the unique solution of
\eqref{rwp1} is given by the function $f$ in the theorem.
\vskip 0.2cm
\noindent
{\bf Step 1.}
Let us check first that $\widetilde A$ is an operator. In fact, if 
\begin{equation*}
 \begin{pmatrix} f \\ k \end{pmatrix}\in\dom\widetilde A\quad\text{and}\quad f=k=0,
\end{equation*} 
then obviously $T^*f=0$ and hence
$\hat f=0$. This yields $\Upsilon_0\hat f=0=\Gamma_0\hat k$ and $\Upsilon_1\hat f=0=\Gamma_1\hat k$. Therefore 
$\hat k=\{0,k^\prime\}\in S$ and as $S$ is an operator $k^\prime=0$ follows. 
The fact that $\widetilde A$ is symmetric in the Krein space $L^2(\partial\Omega)\times \cK$ follows from 
the special form of $\dom \widetilde A$ and the identities \eqref{greenabs} and \eqref{green} for the ordinary 
boundary triple $\{L^2(\partial\Omega),\Upsilon_0,\Upsilon_1\}$ and the generalized boundary triple $\{L^2(\partial\Omega),\Gamma_0,\Gamma_1\}$. Indeed, for $\bigl(\begin{smallmatrix} f \\ k\end{smallmatrix}\bigr),\, \bigl(\begin{smallmatrix} g \\ h\end{smallmatrix}\bigr)\in\dom\widetilde A$ we have $\Upsilon_0\hat f=\Gamma_0\hat k$, 
$\Upsilon_0\hat g=\Gamma_0\hat h$, $\Upsilon_1\hat f=-\Gamma_1\hat k$, $\Upsilon_1\hat g=-\Gamma_1\hat h$ and hence
\begin{equation*}
\begin{split}
&\left[\widetilde A\begin{pmatrix}f\\k\end{pmatrix},\begin{pmatrix}g\\h\end{pmatrix}\right]-
\left[\begin{pmatrix}f\\k\end{pmatrix},\widetilde A\begin{pmatrix}g\\h\end{pmatrix}\right]\\
&\quad
=(\Upsilon_1\hat f,\Upsilon_0\hat g)-(\Upsilon_0\hat f,\Upsilon_1\hat g)
+(\Gamma_1\hat k,\Gamma_0\hat h)-(\Gamma_0\hat k,\Gamma_1\hat h)=0.
\end{split}
\end{equation*}

In order to prove that $\widetilde A$ is selfadjoint in $L^2(\Omega)\times\cK$ it is sufficient to verify
that the operators $\widetilde A-\mu$ and $\widetilde A-\bar\mu$ are surjective for some $\mu\in\cU$. 
We show only $\ran(\widetilde A-\mu)=
L^2(\Omega)\times\cK$, the same reasoning applies for $\widetilde A-\bar\mu$.
By assumption 
$\mu\in\cO$ is such that $(M(\mu)+\tau(\mu))^{-1}\in\cL(L^2(\partial\Omega))$ and moreover, $\mu$ belongs to
$\rho(T_D)\cap\rho(S_0)$ as $\sigma(T_D)\subset\dR$ and $\tau$ is holomorphic on $\cO\cup\cO^*$.
Let $g\in L^2(\Omega)$, $h\in\cK$ and define
$\hat f=\{f,\mu f+g\}$ and $\hat k=\{k ,\mu k+h\}$ by
\begin{equation}\label{f}
f:=(T_D-\mu)^{-1}g-\gamma(\mu)\bigl(M(\mu)+\tau(\mu)\bigr)^{-1}(\gamma(\bar\mu)^*g+\gamma_\tau(\bar\mu)^+h)\in L^2(\Omega)
\end{equation}
and 
\begin{equation*}
k:=(S_0-\mu)^{-1}h-\gamma_\tau(\mu)\bigl(M(\mu)+\tau(\mu)\bigr)^{-1}(\gamma(\bar\mu)^*g+\gamma_\tau(\bar\mu)^+h)\in\cK.
\end{equation*}
Here $\gamma$ is the $\gamma$-field of the ordinary boundary triple $\{L^2(\partial\Omega),\Upsilon_0,\Upsilon_1\}$ and $\gamma_\tau$ is the 
$\gamma$-field corresponding to the generalized boundary triple $\{L^2(\partial\Omega),\Gamma_0,\Gamma_1\}$.
Note that $\hat f\in T^*$ since $\gamma(\mu)(M(\mu)+\tau(\mu))^{-1}(\gamma(\bar\mu)^*g+\gamma_\tau(\bar\mu)^+h)\in\cN_{\mu,T^*}$ and
\begin{equation}\label{ftd}
\bigl\{(T_D-\mu)^{-1}g,(I+\mu(T_D-\mu)^{-1})g\bigr\}\in T_D.
\end{equation}
An analogous argument shows $\hat k\in \dom\Gamma\subset S^+$. We claim that $\{\hat f,\hat k\}$ satisfies the boundary conditions 
$\Upsilon_0\hat f =  \Gamma_0\hat k$ and $\Upsilon_1\hat f = - \Gamma_1\hat k$, so that 
$\bigl(\begin{smallmatrix} f\\ k \end{smallmatrix}\bigr)$ belongs to $\dom\widetilde A$. In fact, 
as $T_D=\ker\Upsilon_0$ it follows from \eqref{f}, \eqref{ftd} and \eqref{gambar} that 
\begin{displaymath}
\begin{split}
 \Upsilon_0\hat f&= -\bigl(M(\mu)+\tau(\mu)\bigr)^{-1}(\gamma(\bar\mu)^*g+\gamma_\tau(\bar\mu)^+h),\\
 \Upsilon_1\hat f&= \gamma(\bar\mu)^*g-M(\mu)\bigl(M(\mu)+\tau(\mu)\bigr)^{-1}(\gamma(\bar\mu)^*g+\gamma_\tau(\bar\mu)^+h),
\end{split}
\end{displaymath}
and analogously,
\begin{displaymath}
\begin{split}
 \Gamma_0\hat k&= -\bigl(M(\mu)+\tau(\mu)\bigr)^{-1}(\gamma(\bar\mu)^*g+\gamma_\tau(\bar\mu)^+h),\\
 \Gamma_1\hat k&= \gamma_\tau(\bar\mu)^+h-\tau(\mu)\bigl(M(\mu)+\tau(\mu)\bigr)^{-1}(\gamma(\bar\mu)^*g+\gamma_\tau(\bar\mu)^+h).
\end{split}
\end{displaymath}
Hence we have $\Upsilon_0\hat f =  \Gamma_0\hat k$ and 
\begin{equation*}
\begin{split}
\Upsilon_1\hat f=& \gamma(\bar\mu)^*g-(\gamma(\bar\mu)^*g+\gamma_\tau(\bar\mu)^+h)\\
&\qquad+\tau(\mu)\bigl(M(\mu)+
\tau(\mu)\bigr)^{-1}(\gamma(\bar\mu)^*g+\gamma_\tau(\bar\mu)^+h)
=-\Gamma_1\hat k,
\end{split}
\end{equation*}
i.e., $\{\hat f,\hat k\}\in\widetilde A$ and it follows that
\begin{equation*}
(\widetilde A-\mu)\begin{pmatrix} f\\k\end{pmatrix}=
\begin{pmatrix} \mu f + g\\ \mu k+h\end{pmatrix}-\mu\begin{pmatrix} f\\k\end{pmatrix}=\begin{pmatrix} g\\h\end{pmatrix}
\end{equation*}
holds. As the elements $g\in L^2(\Omega)$ and $h\in\cK$ were chosen arbitrary we conclude 
$\ran(\widetilde A-\mu)=L^2(\Omega)\times\cK$.

\vskip 0.2cm
\noindent
{\bf Step 2.} Next it will be verified that for $\lambda\in\cU$ the unique solution of \eqref{rwp1} is given by
\begin{equation}\label{fsol}
f=P_{L^2}(\widetilde A-\lambda)^{-1}\begin{pmatrix} g \\ 0\end{pmatrix}.
\end{equation}
We note first that the set $\cU$ is a subset of $\rho(\widetilde A)$. In fact, for every $\lambda\in\cU$ the same argument
as in Step 1 of the proof shows that $\widetilde A-\lambda$ and $\widetilde A-\bar\lambda$ are surjective and hence 
$\ker(\widetilde A-\bar\lambda)=\{0\}=\ker(\widetilde A-\lambda)$, i.e. $\lambda,\bar\lambda\in\rho(\widetilde A)$.
For $f$ in \eqref{fsol} we have
\begin{equation*}
(\widetilde A-\lambda)^{-1}\begin{pmatrix} g \\ 0\end{pmatrix}=\begin{pmatrix} f \\ k\end{pmatrix},\quad\text{where}\quad
k:=P_\cK(\widetilde A-\lambda)^{-1}\begin{pmatrix} g \\ 0\end{pmatrix},
\end{equation*}
and from $\widetilde A\subset T^*\times \dom\Gamma$ and
\begin{equation*}
 \widetilde A\begin{pmatrix} f \\ k\end{pmatrix}=\begin{pmatrix} g \\ 0\end{pmatrix}+\lambda \begin{pmatrix} f \\ k\end{pmatrix}
=\begin{pmatrix} g+\lambda f\\ \lambda k\end{pmatrix} 
\end{equation*}
we conclude that $T^* f=g+\lambda f$ and $k\in\cN_{\lambda,S^+}=\ker(S^+-\lambda)$ holds. 
As $\tau$ is the Weyl function corresponding to 
the generalized boundary triple $\{L^2(\partial\Omega),\Gamma_0,\Gamma_1\}$ and $S^+$ it follows that $\hat k=\{k,\lambda k\}\in\widehat\cN_{\lambda,S^+}\cap\dom\Gamma$ satisfies
$\tau(\lambda)\Gamma_0\hat k=\Gamma_1\hat k$. Therefore, making use of the specific form of $\dom\widetilde A$ and the ordinary boundary triple in
Proposition~\ref{elbt} we obtain
\begin{equation*}
 \tau(\lambda)\iota_- f\vert_{\partial\Omega}=\tau(\lambda)\Upsilon_0\hat f=\tau(\lambda)\Gamma_0\hat k
=\Gamma_1\hat k=-\Upsilon_1\hat f=\iota_+\frac{\partial f_D}{\partial\nu_\ell}\Bigl|_{\partial\Omega}.
\end{equation*}
Hence \eqref{fsol} is a solution of the boundary value problem \eqref{rwp1}. The fact that the compression of the resolvent
of $\widetilde A$ onto $L^2(\Omega)$ has the asserted form follows from Step 1 of the proof by setting $\hat f=\{f,\lambda f+g\}$
and $\hat k=\{k,\lambda k\}$. In this case \eqref{f} reduces to
\begin{equation*}
f=(T_D-\lambda)^{-1}g-\gamma(\lambda)\bigl(M(\lambda)+\tau(\lambda)\bigr)^{-1}\gamma(\bar\lambda)^*g
\end{equation*}
and coincides with $P_{L^2}(\widetilde A-\lambda)^{-1}\vert_{L^2}g$ by \eqref{fsol}.

Finally, we check that for $\lambda\in\cU$ the solution $f$ of \eqref{rwp1} in \eqref{fsol} is unique. Assume that $f_1\in\cD_\max$
is also a solution of \eqref{rwp1}. Then $f-f_1\in\cN_{\lambda,T^*}$ and as $M$ is the Weyl function of $\{L^2(\partial\Omega),\Upsilon_0,\Upsilon_1\}$ we have
\begin{equation*}
M(\lambda)\Upsilon_0(\hat f-\hat f_1)=\Upsilon_1(\hat f-\hat f_1),\qquad \hat f=\{f,T^* f\},\,\hat f_1=\{f_1,T^* f_1\}.
\end{equation*}
On the other hand, since $f$ and $f_1$ both satisfy the boundary condition in \eqref{rwp1} it is clear that 
$\tau(\lambda)\Upsilon_0(\hat f-\hat f_1)=-\Upsilon_1(\hat f-\hat f_1)$ holds and this implies
\begin{equation*}
(M(\lambda)+\tau(\lambda))\Upsilon_0(\hat f-\hat f_1)=0.
\end{equation*}
Since $\lambda\in\cU$ we conclude $\Upsilon_0(\hat f-\hat f_1)=0$, i.e., $\hat f-\hat f_1\in T_D=\ker\Upsilon_0$. 
From $\lambda\in\rho(T_D)$ we then obtain $\hat f=\hat f_1$ and hence the solution $f$ in \eqref{fsol} is unique.
This completes the proof of Theorem~\ref{rwpthm}.
\end{proof}

\begin{remark}
The method applied in the proof of Theorem~\ref{rwpthm} differs from the coupling techniques in \cite[Theorem 4.3]{B08} and \cite[$\S$ 5.2]{DHMS00},
where only ordinary boundary triples were used.
The principal difficulty here is to ensure selfadjointness of $\widetilde A$, a fact that follows immediately via the abstract
boundary condition in \cite{B08,DHMS00}.
\end{remark}

In the special case that $\tau$ in \eqref{rwp1} is a (in general non-strict) $\cL(L^2(\partial\Omega))$-valued Nevanlinna function
the condition $0\in\rho(M(\mu)+\tau(\mu))$ in Theorem~\ref{rwpthm} is automatically satisfied for every nonreal $\mu$, because 
the imaginary part of the Weyl function $M$ of the ordinary boundary triple $\{L^2(\partial\Omega),\Upsilon_0,\Upsilon_1\}$ for $T^*$ 
is uniformly positive (uniformly negative) for $\lambda\in\dC^+$ ($\lambda\in\dC^-$, respectively). This proves the following
corollary.

\begin{corollary}\label{rwpcor}
Assume that the function $\tau$ in the boundary condition in \eqref{rwp1} belongs to the class $N_0(\cL(L^2(\partial\Omega)))$
and let $\{L^2(\partial\Omega),\Upsilon_0,\Upsilon_1\}$ be the ordinary boundary triple for $T^*$ from Proposition~\ref{elbt}
with corresponding $\gamma$-field $\gamma$ and Weyl function $M$.
Then the operator $\widetilde A$ in Theorem~\ref{rwpthm}
is a selfadjoint extension of $T$ in $L^2(\Omega)\times\cK$ and for every $\lambda\in\dC\backslash\dR$ 
the unique solution of 
the boundary value problem \eqref{rwp1} is given by \eqref{sol}.
\end{corollary}

Observe that for $g=0$ in \eqref{rwp1} and $\lambda\in\cU$ the unique solution
of the homogeneous boundary value problem 
\begin{equation}\label{hombvp}
(\ell -\lambda )f=0\qquad\text{and}\qquad \tau(\lambda)\iota_-f\vert_{\partial\Omega}=
\iota_+\frac{\partial f_D}{\partial\nu_\ell}\bigl|_{\partial\Omega}
\end{equation}
is given by $f=P_{L^2}(\widetilde A-\lambda)^{-1}\vert_{L^2}0=0$, cf. Theorem~\ref{rwpthm}. 
The following proposition shows, roughly speaking, that the nontrivial solutions of the homogeneous problem \eqref{hombvp}
are given by the eigenvalues and eigenvectors of the operator $\widetilde A$.

\begin{proposition}
Let the assumptions be as in Theorem~\ref{rwpthm} and let $\widetilde A$ be the selfadjoint operator in $L^2(\Omega)\times\cK$
from the same theorem. Then the following holds.
\begin{enumerate}
 \item [{\rm (i)}] If $\lambda\in\cO\cup\cO^*$ is an eigenvalue of $\widetilde A$ and 
$\bigl(\begin{smallmatrix} f \\ k\end{smallmatrix}\bigr)\in\ker(\widetilde A-\lambda)$ is a corresponding eigenvector, then
$f\in\cD_\max$ is a nontrivial solution of \eqref{hombvp}.
\item [{\rm (ii)}] If $\lambda\in\cO\cup\cO^*$ and $f\in\cD_\max$ is a nontrivial solution of \eqref{hombvp}, then $\lambda$
is an eigenvalue of $\widetilde A$ and 
$\bigl(\begin{smallmatrix} f \\ k\end{smallmatrix}\bigr)\in\ker(\widetilde A-\lambda)$ for some $k\in\cK$.
\end{enumerate}
\end{proposition}

\begin{proof}
(i) Suppose that $\bigl(\begin{smallmatrix} f \\ k\end{smallmatrix}\bigr)\in\dom\widetilde A$ is an eigenvector corresponding to the eigenvalue
$\lambda\in\cO\cup\cO^*$ of $\widetilde A$. Then we have $\ell f=\lambda f$ and since 
$\hat k=\{k,\lambda k\}\in\widehat\cN_{\lambda,S^+}\cap\dom\Gamma$ it
follows from the specific form of $\dom\widetilde A$ and the fact that $\tau$ is the Weyl function of the generalized
boundary triple $\{L^2(\partial\Omega),\Gamma_0,\Gamma_1\}$ that
\begin{equation*}
\tau(\lambda)\iota_-f\vert_{\partial\Omega}=\tau(\lambda)\Upsilon_0\hat f=\tau(\lambda)\Gamma_0\hat k=\Gamma_1\hat k=
-\Upsilon_1\hat f=\iota_+\frac{\partial f_D}{\partial\nu_\ell}\bigl|_{\partial\Omega}
\end{equation*}
holds. Therefore $f\in\cD_\max$ is a solution of the homogeneous boundary value problem \eqref{hombvp}. 
It remains to show 
$f\not= 0$. Assume the contrary. Then $\hat f=\{f,T^* f\}=0$ and it follows from $0=\Upsilon_0\hat f=\Gamma_0\hat k$
that $\hat k=\{k,\lambda k\}$ belongs to $S_0=\ker\Gamma_0$. 
Since $(\cO\cup\cO^*)\subset\rho(S_0)$ (cf. the beginning of Section~\ref{bvpsec}, Theorem~\ref{mainthm} and Remark~\ref{nonminrem})
we conclude $k=0$, a contradiction to $\bigl(\begin{smallmatrix} f \\ k\end{smallmatrix}\bigr)$ being an eigenvector.

(ii) Let $f\in\cD_\max$ be a nontrivial solution of \eqref{hombvp}. Then the boundary condition 
$\tau(\lambda)\Upsilon_0\hat f=-\Upsilon_1\hat f$, $\hat f=\{f,\lambda f\}$, is fulfilled and as $\lambda\in(\cO\cup\cO^*)\subset\rho(S_0)$,
$S_0=\ker\Gamma_0$, we can decompose $\dom \Gamma$ in the form $\dom\Gamma=S_0\,\widehat +\,\widehat\cN_{\lambda,\dom\Gamma}$, 
cf. \eqref{dect}. Since $\{L^2(\partial\Omega),\Gamma_0,\Gamma_1\}$ is a generalized boundary triple for $S^+=\overline{\dom\Gamma}$
the map $\Gamma_0:\dom\Gamma\rightarrow L^2(\partial\Omega)$ is onto and hence there exists 
$\hat k=\{k,\lambda k\}\in\widehat\cN_{\lambda,\dom\Gamma}=\widehat\cN_{\lambda,S^+}\cap\dom\Gamma$ such that 
$\Gamma_0\hat k=\iota_-f\vert_{\partial\Omega}$ holds. Hence we have $\Gamma_0\hat k=\Upsilon_0\hat f$,
$\tau(\lambda)\Gamma_0\hat k=\Gamma_1\hat k$, and therefore
\begin{equation*}
\Upsilon_1\hat f=-\tau(\lambda)\Upsilon_0\hat f=-\tau(\lambda)\Gamma_0\hat k=-\Gamma_1\hat k,
\end{equation*}
i.e., $\bigl(\begin{smallmatrix} f \\ k\end{smallmatrix}\bigr)\in\dom\widetilde A$ is an eigenvector corresponding to the 
eigenvalue $\lambda$ of $\widetilde A$.
\end{proof}

\subsection{An example: A rational Nevanlinna function $\tau$}\label{tauex}

Let $\alpha_i,\beta_i\in\cL(L^2(\partial\Omega))$, $i=1,\dots,m$, be bounded selfadjoint operators in $L^2(\partial\Omega)$ 
and assume that
$\beta_i\geq 0$ holds for all $i=1,\dots,m$ and $0\in\rho(\beta_1)$.
We consider the boundary value problem \eqref{rwp1} with a  function
$\tau$ of the form
\begin{equation}\label{tau}
\tau(\lambda)=\alpha_1+\lambda \beta_1+\sum_{i=2}^m\beta_i^{1/2}(\alpha_i-\lambda)^{-1}\beta_i^{1/2},\qquad
\lambda\in\bigcap_{i=2}^m\rho(\alpha_i).
\end{equation}
Observe that $\tau$ is an $\cL(L^2(\partial\Omega))$-valued Nevanlinna function with the property $0\in\rho(\Imag \tau(\lambda))$ 
for all $\lambda\in\dC\backslash\dR$ and hence $\tau$ is (uniformly) strict. The next theorem, in which
a solution operator $\widetilde A$ of the boundary value problem \eqref{rwp1}, \eqref{tau} is explicitely constructed, 
is essentially a consequence of Theorem~\ref{rwpthm} and an explicit realization of the function \eqref{tau}
as the Weyl function of an ordinary boundary triple in the product space 
\begin{equation*}
L^2(\partial\Omega)^m=L^2(\partial\Omega)\times\, \dots\,\times  L^2(\partial\Omega)\qquad\text{($m$ copies)}.
\end{equation*}
A special case of Theorem~\ref{ratrwp} below was announced in \cite{BGAMM}.
For ordinary second order differential operators in $L^2(I)$, $I\subset\dR$, and scalar rational Nevanlinna functions in the boundary condition 
a solution operator of similar form in $L^2(I)\oplus\dC^m$ as in the next result can be found in \cite{B03}, see also \cite{BBW02-1,BBW02-2}.

\begin{theorem}\label{ratrwp}
Let $\tau$ be a rational $\cL(L^2(\partial\Omega))$-valued Nevanlinna function of the form \eqref{tau}
and let $\gamma$ and $M$ be as in Proposition~\ref{elbt}. Then
\begin{equation*}
\begin{split}
\widetilde A\begin{pmatrix} f \\ k_1 \\ k_2 \\ \vdots \\ k_m\end{pmatrix} &=\begin{pmatrix} \ell f \\ 
k_1^\prime
 \\ \beta_2^{1/2}\beta_1^{-1/2}k_1+\alpha_2k_2 \\ \vdots \\ \beta_m^{1/2}\beta_1^{-1/2}k_1+\alpha_mk_m \end{pmatrix},\qquad
\begin{matrix} f=f_D+f_\eta\in\cD_\max,\\  \\ k_1,\dots,k_m,k_1^\prime\in L^2(\partial\Omega),\end{matrix}\\
\dom\widetilde A&=\left\{\begin{pmatrix} f \\ k_1 \\ \vdots \\ k_m\end{pmatrix}:
\begin{array}{cl}
\iota_- f_\eta\vert_{\partial\Omega}  &=  \beta_1^{-1/2}k_1\\
\iota_+\frac{\partial f_D}{\partial\nu_\ell}\bigl|_{\partial\Omega} &=
\alpha_1\beta_1^{-1/2}k_1+\beta_1^{1/2}k_1^\prime-\sum_{i=2}^{m}\beta_i^{1/2}k_i
\end{array}
\right\},
\end{split}
\end{equation*}
is a selfadjoint operator in the Hilbert space $L^2(\Omega)\times L^2(\partial\Omega)^m$ and for every 
$\lambda$ in $\rho(\widetilde A)\cap\rho(T_D)\cap\mathfrak h(\tau)$ the unique solution of the boundary value problem \eqref{rwp1}
is given by \eqref{sol}.
\end{theorem}

\begin{proof}
The statements in Theorem~\ref{ratrwp} will follow by applying Theorem~\ref{rwpthm} to an explicit 
realization of the function $\tau$ in \eqref{tau} 
as the Weyl function of an ordinary boundary triple $\{L^2(\partial\Omega),\Gamma_0,\Gamma_1\}$ for some closed symmetric operator 
in $L^2(\partial\Omega)^m$.

Denote the functions $k\in L^2(\partial\Omega)^m$ in the form $k=(k_1,\dots,k_m)^\top$, $k_i\in L^2(\partial\Omega)$,
$i=1,\dots,m$, and consider the non-densely defined operator
\begin{equation*}
 \begin{split}
 S(k_1,\dots,k_m)^\top&=\Bigl(\sum_{i=2}^m\beta_1^{-1/2}\beta_i^{1/2} k_i,\alpha_2 k_2,\dots,\alpha_m k_m\Bigr)^{\top},\\
 \dom S&=\bigl\{(k_1,\dots,k_m)^{\top}\in L^2(\partial\Omega)^m : k_1=0\bigr\},
 \end{split}
\end{equation*}
in $L^2(\partial\Omega)^m$. The scalar products in $L^2(\partial\Omega)$ and $L^2(\partial\Omega)^m$ will both be denoted 
by $(\cdot,\cdot)$. We hope that this does not lead to any confusion.
As $\alpha_i=\alpha_i^*$, $i=1,\dots,m$, it follows that $(Sk,k)$ is real for all $k\in\dom S$ and hence $S$ is symmetric. We claim that
the adjoint of $S$ is given by 
\begin{equation}\label{sstern}
S^*=\left\{\left\{\begin{pmatrix}k_1 \\ k_2 \\ \vdots \\ k_m \end{pmatrix},
\begin{pmatrix} k_1^\prime \\ \beta_2^{1/2}\beta_1^{-1/2}k_1+\alpha_2k_2 \\ \vdots \\ \beta_m^{1/2}\beta_1^{-1/2}k_1+\alpha_mk_m
\end{pmatrix}\right\}: k_1,\dots, k_m,k_1^\prime \in L^2(\partial\Omega)\right\}.
\end{equation}
In fact, for $l\in\dom S$ and an element $\hat k=\{k,k^\prime\}$ belonging to the right hand side of \eqref{sstern} we compute
\begin{equation*}
 \begin{split}
  (Sl,k)-(l,k^\prime)&=\sum_{i=2}^m(\beta_1^{-1/2}\beta_i^{1/2}l_i,k_1)+\sum_{i=2}^m(\alpha_il_i,k_i)\\
&\qquad\qquad-
\sum_{i=2}^m(l_i,\beta_i^{1/2}\beta_1^{-1/2}k_1+\alpha_ik_i)=0
 \end{split}
\end{equation*}
and hence the right hand side of \eqref{sstern} is a subset of $S^*$. Furthermore, for each $l\in\dom S$ and $\hat k=\{k,k^\prime\}\in S^*$ we have
\begin{equation*}
0=(Sl,k)-(l,k^\prime)=\sum_{i=2}^m(\beta_1^{-1/2}\beta_i^{1/2}l_i,k_1)+\sum_{i=2}^m(\alpha_il_i,k_i)-\sum_{i=2}^m(l_i,k_i^\prime).
\end{equation*}
Therefore, by inserting $l=(0,\dots,0,l_j,0,\dots,0)^\top$, $l_j\in L^2(\partial\Omega)$, $j=2,\dots,m$, we obtain
\begin{equation*}
k_j^\prime=\beta_j^{1/2}\beta_1^{-1/2}k_1+\alpha_j k_j,\qquad j=2,\dots,m,
\end{equation*}
i.e., $S^*$ is a subset of the right hand side of \eqref{sstern} and hence $S^*$ is given by \eqref{sstern}. 

Let us check that $\{L^2(\partial\Omega),\Gamma_0,\Gamma_1\}$, where
\begin{equation*}
\Gamma_0\hat k=\beta_1^{-1/2}k_1\quad\text{and}\quad
\Gamma_1\hat k=\alpha_1\beta_1^{-1/2}k_1+\beta_1^{1/2}k_1^\prime-\sum_{i=2}^{m}\beta_i^{1/2}k_i,\quad\hat k\in S^*,
\end{equation*}
is an ordinary boundary triple for $S^*$ with $\tau$ in \eqref{tau} as corresponding Weyl function. Since for an element 
$\hat k=\{k,k^\prime\}\in S^*$ the entries
$k_1$ and $k_1^\prime$ are arbitrary elements in $L^2(\partial\Omega)$ it follows immediately from $0\in\rho(\beta_1)$
that the mapping $(\Gamma_0,\Gamma_1)^\top:S^*\rightarrow L^2(\partial\Omega)\times L^2(\partial\Omega)$ is onto. Next we verify the identity \eqref{greenexp}. For 
$\hat l=\{l,l^\prime\}$ and $\hat k=\{k,k^\prime\}\in S^*$ a straightforward computation shows
\begin{equation*}
\begin{split}
 (l^\prime,k)&-(l,k^\prime)=(\beta_1^{1/2}l_1^\prime,\beta_1^{-1/2}k_1)-(\beta_1^{-1/2}l_1,\beta_1^{1/2}k_1^\prime)\\
&\qquad\qquad+\sum_{i=2}^m(\beta_i^{1/2}\beta_1^{-1/2}l_1+\alpha_il_i,k_i) - \sum_{i=2}^m(l_i,\beta_i^{1/2}\beta_1^{-1/2}k_1+\alpha_ik_i)\\
&=\left(\beta_1^{1/2}l_1^\prime-\sum_{i=2}^m\beta_i^{1/2}l_i\,,\,\beta_1^{-1/2}k_1\right)-
\left(\beta_1^{-1/2}l_1\,,\,\beta_1^{1/2}k_1^\prime-\sum_{i=2}^m\beta_i^{1/2}k_i\right)\\
&=(\Gamma_1\hat l,\Gamma_0\hat k)-(\Gamma_0\hat l,\Gamma_1\hat k),
\end{split}
\end{equation*}
where we have used $\alpha_1=\alpha_1^*$ in the last step.
Observe that the selfadjoint relation $S_0=\ker\Gamma_0$ is given by
\begin{equation*}
S_0=\Bigl\{\bigl\{(0,k_2\dots,k_m)^\top,(k_1^\prime,\alpha_2k_2,\dots,\alpha_m k_m)^\top\bigr\}:
k_1^\prime,k_2,\dots,k_m\in L^2(\partial\Omega)\Bigr\}
\end{equation*}
and that for $\lambda\in\rho(S_0)=\bigcap_{i=2}^m \rho(\alpha_i)$ the resolvent of $S_0$ is a diagonal block operator
matrix in $L^2(\partial\Omega)^m$ with entries $0,(\alpha_2-\lambda)^{-1},\dots,(\alpha_m-\lambda)^{-1}$ on the diagonal.
Let now $\hat k=\{k,\lambda k\}\in\widehat\cN_{\lambda,S^*}$ and $\lambda\in\rho(S_0)$. Then we have
\begin{equation*}
k_1^\prime=\lambda k_1\qquad\text{and}\qquad\beta_i^{1/2}\beta_1^{-1/2}k_1=(\lambda-\alpha_i)k_i,\,\,\,\,i=2,\dots,m,
\end{equation*}
and this implies
\begin{equation*}
\begin{split}
&\left(\alpha_1+\lambda \beta_1+\sum_{i=2}^m\beta_i^{1/2}(\alpha_i-\lambda)^{-1}\beta_i^{1/2}\right)\Gamma_0\hat k\\
&\qquad\qquad\qquad
=\alpha_1\beta_1^{-1/2}k_1+\lambda\beta_1^{1/2}k_1+\sum_{i=2}^m\beta_i^{1/2}(\alpha_i-\lambda)^{-1}\beta_i^{1/2}\beta_1^{-1/2}k_1\\
&\qquad\qquad\qquad=\alpha_1\beta_1^{-1/2}k_1+\beta_1^{1/2}k_1^\prime-\sum_{i=2}^m\beta_i^{1/2}k_i=\Gamma_1\hat k
\end{split}
\end{equation*}
for $\lambda\in\rho(S_0)$. Hence $\tau$ is the Weyl function of the ordinary 
boundary triple $\{L^2(\partial\Omega),\Gamma_0,\Gamma_1\}$.

Now we apply Theorem~\ref{rwpthm} to the present situation. It follows directly from \eqref{sstern} and the 
definition of the boundary triples $\{L^2(\partial\Omega),\Upsilon_0,\Upsilon_1\}$ in Proposition~\ref{elbt} and
$\{L^2(\partial\Omega),\Gamma_0,\Gamma_1\}$ above that the solution operator $\widetilde A$ in Theorem~\ref{rwpthm}
has the asserted form. As $\tau$ is a Nevanlinna function $\dC\backslash\dR$ is subset of $\cU$, cf. the consideration before
Corollary~\ref{rwpcor}, and hence for every $\lambda\in\dC\backslash\dR$ the unique solution $f\in\cD_\max$ 
of \eqref{rwp1} is given by \eqref{sol}. It can be shown with similar arguments as in step 1 of the 
proof of Theorem~\ref{rwpthm} that this is also true on the larger set $\rho(\widetilde A)\cap\rho(T_D)\cap\mathfrak h(\tau)$.
\end{proof}

In the next corollary we consider the special case of a linear $\cL(L^2(\partial\Omega))$-valued Nevanlinna function $\tau$
in the boundary condition of \eqref{rwp1}.
Similar
$\lambda$-linear elliptic boundary value problems were investigated in, e.g., \cite{BHLN01,ES65-1,ES65-2}.

\begin{corollary}\label{lincor}
Let $\alpha,\beta$ be bounded selfadjoint operators in $L^2(\partial\Omega)$ and assume that $\beta$  is uniformly positive. 
Then
\begin{equation*}
\begin{split}
\widetilde A\begin{pmatrix} f \\ k\end{pmatrix} &=\begin{pmatrix} \ell f \\ 
\beta^{-1/2}\iota_+\frac{\partial f_D}{\partial\nu_\ell}\vert_{\partial\Omega} -\beta^{-1/2}\alpha\beta^{-1/2}k
\end{pmatrix}\\
\dom\widetilde A&=\left\{\begin{pmatrix} f \\ k\end{pmatrix}\in\cD_\max\times L^2(\partial\Omega):
\iota_- f_\eta\vert_{\partial\Omega}=\beta^{-1/2}k\right\}
\end{split}
\end{equation*}
is a selfadjoint operator in $L^2(\Omega)\times L^2(\partial\Omega)$ and for $g\in L^2(\Omega)$ and
$\lambda\in\rho(\widetilde A)\cap\rho(T_D)$ 
the unique solution $f\in\cD_\max$ of the boundary value problem 
\begin{equation*}
(\ell -\lambda)f=g,\qquad (\alpha+\lambda\beta) \iota_- f\vert_{\partial\Omega}=\iota_+\frac{\partial f_D}{\partial\nu_\ell}
\bigl|_{\partial\Omega},
\end{equation*}
is given by \eqref{sol}.
\end{corollary}

\end{document}